\documentclass[12pt]{article}
\usepackage{amsmath}
\usepackage{graphicx}
\usepackage{enumerate}
\usepackage{natbib}
\usepackage{url} 
\usepackage{subfigure}
\usepackage[figuresright]{rotating}
\usepackage{bm}
\usepackage{float}
\usepackage{xcolor}
\RequirePackage{amsthm,amsmath,amsfonts,amssymb}
\theoremstyle{plain}
\newtheorem{thm}{Theorem}[section]
\newtheorem{remark}{Remark}[section]

\numberwithin{equation}{section}
\newcommand{\blind}{0}

\begin{document}

\def\spacingset#1{\renewcommand{\baselinestretch}%
{#1}\small\normalsize} \spacingset{1}

\if0\blind
{
  \title{\bf Average quantile regression: a new non-mean regression model and coherent risk measure}
   \author{Rong Jiang \hspace{.2cm}\\
    Shanghai University of International Business and Economics, China\\
    and \\
    M.C. Jones \\
    The Open University, UK\\
     and \\
    Keming Yu \\
    Brunel University, UK\\
     and \\
    Jiangfeng Wang \\
    Zhejiang Gongshang University, China}
  \maketitle
} \fi

\if1\blind
{
  \bigskip
  \bigskip
  \bigskip
  \begin{center}
    {\LARGE\bf Title}
\end{center}
  \medskip
} \fi

\bigskip
\begin{abstract}
Regression models that go beyond the mean, alongside coherent risk measures, have been important tools in modern data analysis. This paper introduces the innovative concept of Average Quantile Regression (AQR), which is smooth at the quantile-like level, comonotonically additive, and explicitly accounts for the severity of tail losses relative to quantile regression. AQR serves as a versatile regression model capable of describing distributional information across all positions, akin to quantile regression, yet offering enhanced interpretability compared to expectiles. Numerous traditional regression models and coherent risk measures can be regarded as special cases of AQR. As a flexible non-parametric regression model, AQR demonstrates outstanding performance in analyzing high-dimensional and large datasets, particularly those generated by distributed systems, and provides a convenient framework for their statistical analysis. The corresponding estimators are rigorously derived, and their asymptotic properties are thoroughly developed. In a risk management context, the case study confirms AQR's effectiveness in risk assessment and portfolio optimization.
\end{abstract}

\noindent%
\vfill

\newpage
\spacingset{1.9} 

\section{Introduction}
   Regression analysis is one of the most vital tools in statistical data analysis and plays a significant role across various fields.
   Generalized Linear Models (GLMs), formulated by \cite{r10}, provide a framework that unifies numerous mean based regression models. They extend traditional linear regression by using a link function to connect a linear predictor to the mean of the response variable.
  However, in many disciplinary fields, such as meteorology, life sciences, and financial risk management, the involved data often exhibit skewed distribution and heterogeneity characteristics.
  Regression models that account for heavy tails, asymmetry, and outliers have therefore been attracting significant attention.
  In view of this, in related studies, in addition to focusing on the mean of the conditional distribution of the response given the covariates, researchers are often interested in the tail features of the data.
   \par
   Popular regression models that go beyond the mean estimation include quantile regression \citep{r3} and expectile regression \citep{r43}. Quantile regression (QR) is predominant in the literature due to its excellent interpretability.
   However, the standard quantile regression objective function lacks smoothness, even when the underlying quantile function is absolutely continuous. The absence of a second derivative in the objective function complicates statistical inference and may lead to non-unique solutions. Furthermore, the asymptotic normality of the standard quantile regression estimator relies on Bahadur-Kiefer representations, which are characterized by slow convergence rates, as noted by \cite{r51}. In contrast, the loss function of expectile regression is differentiable everywhere. Additionally, quantile-based Value at Risk (VaR) fails to satisfy the subadditivity property, making it inconsistent as a risk measure according to the axiomatic framework proposed by \cite{r15}. Expectile-based VaR is a coherent risk measure. However, its lack of comonotonic additivity poses significant challenges for regulatory risk standards \citep{r50}. Moreover, its inferences are more sensitive to extreme values or outliers and lack intuitive explanations. Other widely used coherent risk measures, such as Expected Shortfall \citep{r28} and Extremiles \citep{r1}, often exhibit excessive conservatism, limiting their practical applicability for individual financial institutions. While Expected Shortfall emphasizes tail-risk scenarios but lacks the distributional modeling capabilities inherent in regression frameworks. Conversely, Extremiles face challenges in providing clear quantile-level interpretations, despite their theoretical focus on extreme events.
\par
Our work introduces a novel family of regression models that extends beyond mean regression, alongside coherent risk measures. These models are particularly significant in real-world applications where the effects of explanatory variables differ across various levels of outcomes (e.g., income, health, risk) and extreme cases, where traditional mean regression methods, such as ordinary least squares (OLS), prove insufficient. Coherent risk measures, on the other hand, play a crucial role in finance and risk management by offering mathematically rigorous and economically meaningful approaches to assessing financial risks. They address key limitations of older measures like Value-at-Risk (VaR) and ensure consistency in decision-making processes. The proposed models adhere to principles of comonotonic additivity and coherence as law-invariant spectral risk measures, referred to as Average Quantile-like Regression (AQR). Notably, many traditional regression models and risk measures can be viewed as special cases within the AQR framework. Furthermore, this family of models addresses inherent limitations in existing regression techniques and risk measurement approaches. The AQR framework significantly enhances risk measurement methodologies by encompassing a broad spectrum of risk assessments, ranging from measures more conservative than Expected Shortfall to those more aggressive than Extremiles. Empirical case studies demonstrate the superior performance of AQR in portfolio optimization and environmental applications. Additionally, as a flexible family of non-parametric regression models, AQR exhibits remarkable capability in handling large-scale and high-dimensional datasets, especially those generated by distributed systems, while providing a robust framework for statistical analysis.

 \subsection{Average quantile-like regression models}
We extend classical QR  to incorporate newly proposed regression models by utilizing a density function, $\text{J}_\tau(s)$, $0\leq s\leq1$, over which the conditional quantile function is weighted averaged. We refer to the result as AQR. For $0<\tau<1$, the conditional $\tau$-th AQR is defined as follows:
\begin{equation}
	\begin{split}
		\xi_{\tau}(\bm{Y}|\bm{x})=\int_0^1Q_{\bm{Y}|\bm{x}}(s)\text{J}_{\tau}(s)ds,
	\end{split}
\end{equation}
where $\bm{Y}$ is a dependent variable in $\mathbb{R}$, $\bm{X}$ is a vector of covariates in $\mathbb{R}^p$ and $Q_{\bm{Y}|\bm{x}}(s)$ is the quantile function of the conditional distribution of $\bm{Y}$ given $\bm{X}=\bm{x}$. The weight function $\text{J}_{\tau}(s)$ in model (1.1) must satisfy certain conditions labelled as {\bf C1} in Section 2.1.
\par
Notice that the weighting function proposed in this paper incorporates two parameters,
$s$ and $\tau$. The parameter $s$ is used to weight the quantile function, while the parameter $\tau$ allows for fitting different positions of the data, resulting in the average quantile being a further quantile-like function of $\tau$.
For risk measurement, this contrasts with a couple of existing risk measures \citep{r16,r31,r7}, which    focus solely on weighting the quantile function using $s$ for, implicitly, a single value of $\tau$.
 \subsubsection{Regression model examples}
 The classical QR and many quantile-related regression models shown below are special cases of AQR.
\par
(1) Quantile regression:
 $Q_{\bm{Y}|\bm{x}}(\tau)$ is equal to $\xi_{\tau}(\bm{Y}|\bm{x})$ with $\text{J}_{\tau}(s)=\delta(s-\tau)$, where $\delta(\cdot)$ is the Dirac delta function ($\delta(u)=0$ with $u\neq0$ and $\int_{-\infty}^{+\infty}\delta(u)du=1$).
\par
 (2) Extremile regression \citep{r4} given by
 $\xi_{\tau}(\bm{Y}|\bm{x})$ with
   \begin{equation}
	\begin{split}
		\text{J}_{\tau}(s)=\left \{
		\begin{array}{ll}
			r_1(\tau)(1-s)^{r_1(\tau)-1},&\textrm{if}~0< \tau\leq 1/2,\\
			r_2(\tau)s^{r_2(\tau)-1},&\textrm{if}~1/2<\tau<1,
		\end{array}
		\right.
	\end{split}
\end{equation}
  where $r_1(\tau)=r_2(1-\tau)=\log(1/2)/\log(1-\tau)$.
\subsubsection{Risk measure examples}
We further show that $\omega_{\tau}\xi_{\tau}(\bm{Y}|\bm{x})$ based on the proposed AQR includes several of the  most popular risk measures in finance, where
\begin{equation*}
	 \begin{split}
		 \omega_{\tau}=\left \{
		\begin{array}{rl}
			-1,~~&\textrm{if}~0< \tau\leq 1/2,\\
		1,~~&\textrm{if}~1/2<\tau<1.
		\end{array}
		\right.
	\end{split}
\end{equation*}
Some existing special cases of $\omega_{\tau}\xi_{\tau}(\bm{Y}|\bm{x})$ are shown below.
\par
(1) Expected shortfall (ES, \cite{r28})
is equal to $\omega_{\tau}\xi_{\tau}(\bm{Y}|\bm{x})$ with
\begin{equation*}
	\begin{split}
		\text{J}_{\tau}(s)=\left \{
		\begin{array}{ll}
			 \text{I}(s<\tau)/\tau,~~&\textrm{if}~0< \tau\leq 1/2,\\
			\text{I}(s\geq\tau)/(1-\tau),~~&\textrm{if}~1/2<\tau<1.
		\end{array}
		\right.
	\end{split}
\end{equation*}
Here, $\text{I}(\cdot)$ is the indicator function.
\par
(2)
Exponential spectral risk measure \citep{r17} is equal to $\omega_{\tau}\xi_{\tau}(\bm{Y}|\bm{x})$ with
\begin{equation*}
	\begin{split}
		\text{J}_{\tau}(s)=\left \{
		\begin{array}{ll}
			(2\tau)^{s}\log(2\tau)/(2\tau-1),~~&\textrm{if}~0< \tau\leq 1/2,\\
			(2-2\tau)^{1-s}\log(2-2\tau)/(1-2\tau),~~&\textrm{if}~1/2<\tau<1.
		\end{array}
		\right.
	\end{split}
\end{equation*}
\par
(3) The quantile-Value at Risk with $\text{J}_{\tau}(s)=\delta(s-\tau)$ and the signed version of extremile regression with $\text{J}_\tau(s)$ given by (1.2).
\par
Moreover, we will prove in Section 2 that under appropriate $\text{J}_{\tau}(s)$ in condition {\bf C1}, which include those above, $\omega_{\tau}\xi_{\tau}(\bm{Y}|\bm{x})$ is a coherent risk measure \citep{r15}.
A case study in Section 5.1 shows that portfolios based on AQR as a risk measure perform better than existing methods.

\subsection{Nonparametric analysis of distributed data by AQR}
As an application of AQR in modern data analysis, we focus on AQR for distributed systems characterized by large sample sizes and high dimensions. This area has garnered significant research attention in the context of GLMs \citep{r2}, but there has been relatively little exploration of nonparametric models. With regard to the latter, \cite{r21} utilized the one-shot method to investigate single-index composite quantile regression for massive data, while \cite{r27} considered distributed heterogeneous learning based on least squares estimation for generalized partially linear spatially varying coefficient models.
In this paper, we propose a broader class of regression models, distinguishing it from the methods used by \cite{r21} and \cite{r27}.
Subsection 5.3 provides an application of  the proposed distributed data AQR  to a Beijing multi-site air quality dataset and makes a comparison with
an existing nonparametric analysis  of the data by \cite{r34}.

\subsection{The structure of the paper}
Sections 1.1 and 1.2 above outline the main contributions of this paper, with further details provided in the subsequent sections. The proposed AQR generalizes many existing non-mean regression models as well as risk measures, via consideration of the weight function $\text{J}_{\tau}(s)$.  New regression models and risk measures can be proposed within this framework, with examples provided in Section 2. Section 3 focuses on estimation methods for both full data and massive datasets generated by distributed systems.
Section 4 presents simulation examples and demonstrates the application of real data to illustrate the proposed
methods. Finally, we conclude this paper with a brief discussion in Section 5. All technical proofs and an algorithm are included in the Supplementary Material.

\section{New regressions and risk measures from AQR}
In this section, we will explain in detail and propose some new regressions and risk measurement tools. The quantity
 $\omega_{\tau}\xi_{\tau}(\bm{Y}|\bm{x})=\omega_{\tau}\int_0^1Q_{\bm{Y}|\bm{x}}(s)
 \text{J}_{\tau}(s)ds$ is a coherent risk measure \citep{r15} as per Definition A.1 in Appendix A of the Supplementary Material under appropriate choice of $\text{J}_{\tau}(s)$ satisfying condition {\bf C1} below. As mentioned in Section 1.1.2, several classic risk measures are its special cases, such as expected shortfall, exponential spectral risk measure and extremile.

\subsection{Basic properties of AQR}
In order to better understand the basic properties of AQR, note that $F(Q_{\bm{Y}|\bm{x}}(s)|\bm{x})=s$, where $F(\cdot|\bm{x})$ is the conditional distribution of $\bm{Y}$ given $\bm{X}=\bm{x}$. Set $y=Q_{\bm{Y}|\bm{x}}(s)$, then $s=F(y|\bm{x})$. We can rewrite model (1.1) in multiple ways as
\begin{equation}
	\begin{split}
		\xi_{\tau}(\bm{Y}|\bm{x})
=&\int_{-\infty}^{+\infty}
		yf(y|\bm{x})\text{J}_{\tau}\{F(y|\bm{x})\}\,dy=
\text{E}[\bm{Y}\text{J}_{\tau}\{F(\bm{Y}|\bm{x})\}]=\text{E}(\bm{Z}_{\tau}^{\bm{x}}),\\
=&\int_{0}^{+\infty}\left[1-\text{G}_{\tau}\{F(y|\bm{x})\}\right]dy
		-\int_{-\infty}^{0}\text{G}_{\tau}\{F(y|\bm{x})\}dy.
	\end{split}
\end{equation}
Here, $\text{G}_{\tau}(u)=\int_{0}^{u}\text{J}_{\tau}(s)ds$ and the random variable $\bm{Z}_{\tau}^{\bm{x}}$ has cumulative distribution function $\text{F}_{\bm{Z}_{\tau}^{\bm{x}}}(\cdot|\bm{x})=\text{G}_{\tau}\{F(\cdot|\bm{x})\}$.
\par
To establish the basic properties of $\xi_{\tau}(\bm{Y}|\bm{x})$ in (2.1), the following
technical conditions for weighting function $\text{J}_{\tau}(s)$ are needed.
\par
{\bf C1.} $\text{J}_{\tau}(s)$ is (i) positive and bounded ($0\leq\text{J}_{\tau}(s)<\infty$) for all $\tau\in (0,1)$, $s\in [0,1]$ and normalized ($\int_0^1 \text{J}_{\tau}(s)ds=1$); (ii) reverse ($\text{J}_{\tau}(s)=\text{J}_{1-\tau}(1-s)$) and monotonic with respect to $s$ (non-increasing for $0<\tau\leq 1/2$ and non-decreasing for $1/2\leq\tau<1$); (iii) $\text{G}_{\tau}(u)$ is a non-increasing function of $\tau$ for all $u \in (0,1)$.
\begin{remark}
Condition (i) of the function $\text{J}_{\tau}(s)$ is commonly used in weight functions.
$\text{J}_{\tau}(s)$ should be positive because if it displays negative values at some $s$, the corresponding $\omega_{\tau}\xi_{\tau}(\bm{Y}|\bm{x})$ prefers losses to profits at quantile $Q_{\bm{Y}|\bm{x}}(s)$ and is therefore not risk-averse. It can be shown that $\omega_{\tau}\xi_{\tau}(\bm{Y}|\bm{x})$ violates in this case the monotonicity axiom in Definition A.1 in the Appendix. The boundedness condition in (i) guarantees the existence of $\xi_{\tau}(\bm{Y}|\bm{x})$.
Were the normalization condition in (i)
not true, one can show that $\omega_{\tau}\xi_{\tau}(\bm{Y}|\bm{x})$ would violate the translation invariance axiom in Definition A.1. Condition (ii) is used to guarantee the coherence property. ``Reverse'' in condition (ii) is commonly present in risk measurement such as ES and Extremile. The monotonicity in condition (ii) is related to the subadditivity coherency axiom. Condition (iii) ensures that $\xi_{\tau}(\bm{Y}|\bm{x})$ is non-decreasing with respect to $\tau$.
\end{remark}

 \begin{thm}
 Let $\bm{Y}$ given $\bm{X}=\bm{x}$ have a finite absolute first moment and $\text{J}_{\tau}(s)$ satisfy condition {\bf C1}. Then, for any $\tau\in(0,1)$, we have
 \par
 (i) $\xi_{\tau}(\bm{Y}|\bm{x})$ exists and is a non-decreasing function with respect to $\tau$.
Furthermore, if there is an interval $I_u\subset (0,1)$ such that $\text{G}_{\tau}(u)$ is a decreasing function
of $\tau$ for all $u\in I_u$, then $\xi_{\tau}(\bm{Y}|\bm{x})$ increases monotonically with respect to $\tau$, thus avoiding unreasonable crossing that can occur in QR (large quantile is smaller than small quantile).
  \par
 (ii) If the conditional distribution of $\bm{Y}$ given $\bm{X}=\bm{x}$ is symmetric, then
 $\xi_{\tau}(\bm{Y}|\bm{x})-\text{E}(\bm{Y}|\bm{X}=\bm{x})=\text{E}(\bm{Y}|\bm{X}=\bm{x})-\xi_{1-\tau}(\bm{Y}|\bm{x})$, which means that the lower and upper AQR curves are symmetric about the regression mean.
 Furthermore, if $\text{J}_{\tau}(s)=1$ for all $s\in[0,1]$, then $\xi_{\tau}(\bm{Y}|\bm{x})=\text{E}(\bm{Y}|\bm{X}=\bm{x})$. As shown in section 2.2, usually $\text{J}_{0.5}(s)=1$. So under symmetry, such AQR curves are symmetric about the case $\tau=0.5$.
  \par
 (iii) If $\bm{Y}=m(\bm{X})+\sigma(\bm{X})\bm{\varepsilon}$, where $m(\cdot)$ and $\sigma(\cdot)>0$ are unknown functions, and error random variable $\bm{\varepsilon}$ has a finite absolute first moment, then $\xi_{\tau}(\bm{Y}|\bm{x})=m(\bm{x})+\sigma(\bm{x})\xi_{\tau}(\bm{\varepsilon})$, where $\xi_{\tau}(\bm{\varepsilon})=\int_0^1Q_{\bm{\varepsilon}}(s)\text{J}_{\tau}(s)ds$ and $Q_{\bm{\varepsilon}}(s)$ is the quantile function of $\bm{\varepsilon}$.
 This result demonstrates parallel AQR curves under response homogeneity.
  \par
 (iv) $\omega_{\tau}\xi_{\tau}(\bm{Y}|\bm{x})$ is a comonotonically additive coherent risk measure.

\end{thm}
\subsection{New regression and coherent risk measure examples}
In this subsection, we will propose some new coherent risk measures based on equation (1.1), which are also new regression models for fitting tail features. Because of $\text{J}_{\tau}(s)=\text{J}_{1-\tau}(1-s)$ in condition ${\bf C1(ii)}$, we only consider the case of $\tau\in(0,1/2]$, and the case of $\tau\in(1/2,1)$ can be similarly derived.
\subsubsection{Generalized ES}
Note that $\text{J}_{\tau}(s)=\text{I}(0<s<\tau)/\tau$ in $\text{ES}_{\tau}$, which is the density function of Uniform$(0,\tau)$ that provides a constant weight. However, in risk management, the greater the loss, the more attention is paid to it, so it should be given more weight. Therefore, we propose a generalized ES (GES), which gives greater weight to tails away from $\tau$ (see Figure 1a) as follows:
   \begin{equation}
	\begin{split}
		\text{J}_{\tau}(s)=(1+a)\tau^{-1-a}(\tau-s)^a\times\text{I}(0<s<\tau),
	\end{split}
\end{equation}
where $a\geq0$ is a constant. When $a=0$, we have $\text{ES}_{\tau}$.
Moreover, when $a=1$, $\text{J}_{\tau}(s)=2\tau^{-2}(\tau-s)
\times\text{I}(0<s<\tau)$, which is a decreasing linear density on $(0,\tau)$, while
$\text{J}_{\tau}(s)$ with $a=2$, that is,
$3\tau^{-3}(\tau-s)^2\times\text{I}(0<s<\tau)$, is a decreasing quadratic density on $(0,\tau)$. In fact, this ${\text J}_\tau(s)$ in (2.2) is the density of $\tau$ times a ${\rm Beta}(1,a+1)$ random variable.
 Its distribution function is ${\text G}_\tau(u) = 1-(1-u/\tau)^{a+1}\times\text{I}(0<u<\tau)$ which is clearly a decreasing function of $\tau$.
In equation (2.2), $\tau$ is analogous to the quantile level in quantile regression.
Accordingly, $\xi_{\tau}(\bm{Y}|\bm{x})$ is a weighted average of values under the $\tau$ conditional quantile.
\par
As $a$ increases, $\text{J}_{\tau}(s)$ has a larger value in the tail (see the left side of
Figure 1a), which also results in a larger value for $\omega_{\tau}\xi_{\tau}(\bm{Y}|\bm{x})$ (see Theorem 2.2(i) and the right side of Figure 1a), where GES1 and GES2 are GES with $a=1$ and 2.
Moreover, the values of  $\omega_{\tau}\xi_{\tau}(\bm{Y}|\bm{x})$ with GES ($a=0,1,2$) are all larger than QR. Indeed, this is true of any AQR where, for $0<\tau<1/2$, $\text{J}_\tau(s)$ has  support $(0,\tau)$: for such $\text{J}_\tau(s)$, $\xi_{\tau} (\bm{Y}|\bm{x})
=\int_0^\tau Q_{\bm{Y}|\bm{x}}(s)\text{J}_{\tau}(s)ds \leq Q_{\bm{Y}|\bm{x}}(\tau) \int_0^\tau \text{J}_{\tau}(s)ds \leq Q_{\bm{Y}|\bm{x}}(\tau)$ (and similarly for $\tau>1/2$). Therefore, users can choose the appropriate $a$ based on their risk preferences. A higher value means a more cautious approach.

\subsubsection{Generalized Extremile}
Note that, for $0<\tau<1/2$, $\text{J}_{\tau}(s)=r_1(\tau)(1-s)^{r_1(\tau)-1}$ where $r_1(\tau) = \log(1/2)/\log(1-\tau)$ in Extremile, which is the density function of the Beta$(1,r_1(\tau))$ distribution. Therefore, we propose a generalized Extremile (GE) as follows:
\begin{equation}
	\begin{split}
\text{J}_{\tau}(s)=(1+\alpha_{\tau})(1-s)^{\alpha_{\tau}}.
 	\end{split}
\end{equation}
$\text{J}_{\tau}(s)$ is the density function of Beta$(1,\alpha_{\tau}+1)$.
The corresponding distribution function is ${\text G}_\tau(u) = 1-(1-u)^{\alpha_\tau+1}$ which is readily seen to be increasing in $\alpha_\tau$. Therefore, $\alpha_{\tau}: (0,1/2)\rightarrow (0,+\infty)$ should be a decreasing function of $\tau$ to meet condition {\bf C1}(iii). With this choice of $\alpha_\tau$, when $\tau=0.5$, $\text{J}_{0.5}(s)=1$ for all $s\in(0,1)$ so that
$\xi_{\tau}(\bm{Y}|\bm{x})=\text{E}(\bm{Y}|\bm{X}=\bm{x})$.
Extremile is a special case with $\alpha_{\tau}=-\log(2-2\tau)/\log(1-\tau)$.
\par
 When $\alpha_{\tau}$ is an integer, we can obtain that
$ \xi_{\tau}(\bm{Y}|\bm{x})=\text{E}\left\{\min(\bm{Y}_{\bm{x}}^1,\ldots,\bm{Y}_{\bm{x}}^{1+\alpha_{\tau}})\right\},
$
where $\bm{Y}_{\bm{x}}^i$ is the $i$th sample drawn from the conditional distribution of $\bm{Y}$ given ${\bm{X}}={\bm{x}}$.
For the case where $\alpha_{\tau}$ is not an integer, $\xi_{\tau}(\bm{Y}|\bm{x})$ is between the expectations of the minimum of $1+[\alpha_{\tau}]$ and $2+[\alpha_{\tau}]$ independent copies of $\bm{Y}_{\bm{x}}$, where $[\cdot]$ denotes the integer part.
The choice $\alpha_{\tau}=0.5\tau^{-1}-1$ in equation (2.3) is attractive because
 the role of $\tau$ in $\xi_{\tau}(\bm{Y}|\bm{x})$ is then to take $0.5\tau^{-1}$ independent copies of $\bm{Y}_{\bm{x}}$ (see Table 1 to follow).  And Figure 2 to follow shows the influence of different $\tau$ on $\text{J}_{\tau}(\cdot)$.
\par
From Figure 1b, the tails of $\text{J}_{\tau}(s)$ based GEs (GE1 and GE2 are GE with $\alpha_{\tau}=0.5\tau^{-1}-1$ and another alternative choice $0.5\pi\cot(\pi\tau)$) are smaller than that of Extremile, and therefore $\omega_{\tau}\xi_{\tau}(\bm{Y}|\bm{x})$ based on Extremile for the normal distribution is larger than GEs. And interestingly, Extremile is greater than QR, while GEs are less than QR.
Naturally, $\xi_{\tau}(\bm{Y}|\bm{x})$ with GE, like Extremile regression, can also serve as a  regression model.

\subsubsection{Coherent risk measure with truncated Cauchy density function}
The Cauchy distribution is a common and important distribution. We construct a new coherent risk measure with truncated density of the Cauchy
$(0,\alpha_{\tau}^{-1})$ distribution over $(0,1)$ with parameter $\alpha_{\tau}^{-1}$ as follows:
\begin{equation}
	\begin{split}
\text{J}_{\tau}(s)=\frac{\alpha_{\tau}^{-1}}{\alpha_{\tau}^{-2}+s^2}\times\frac{1}{\arctan(\alpha_{\tau})}.
	\end{split}
\end{equation}
In this case,  ${\text G}_\tau(s) = f(\alpha_\tau s)/f(\alpha_\tau)$, where $f(t)=\arctan(t)$ and $t>0$. ${\text G}_\tau(s)$ is increasing with respect to $\alpha_\tau$, because $tf'(t)/f(t)=t/\{(1+t^2)\arctan(t)\}$ can be shown to be decreasing
in $t$.
 When $\tau=0.5$, the $\xi_{\tau}(\bm{Y}|\bm{x})$ with $\text{J}_{\tau}(s)$ in (2.4) is, again, the conditional expectation of $\bm{Y}$ given $\bm{x}$.
In (2.4),
$\alpha_{\tau}^{-1}$ is, of course, a scaling parameter, specifically equal to half the width at half the maximum value of the density.
So $\tau$ is the parameter that controls the scale. For instance, $\tau=0.5(1+1/\alpha_{\tau}^{-1})^{-1}$ under $\alpha_{\tau}=0.5\tau^{-1}-1$. Therefore, the value of $\tau$ can be matched up with the scale parameter $\alpha_{\tau}^{-1}$ (see Table 1).
 Furthermore, the smaller the $\tau$, the greater the weight of the tail (see Figure 2).
 Although in finance, the Cauchy distribution is often used to simulate tail risk due to its heavy tail, a form like the new risk measurement proposed in this article based on equation (2.4) is rare.
\par
The $\text{J}_{\tau}(s)$ in (2.4) and the risk measure $\xi_{\tau}(\bm{Y}|\bm{x})$ with truncated Cauchy distribution (TCRM) based on $\text{J}_{\tau}(s)$ are shown in Figure 1c (TCRM1, TCRM2 and TCRM3 are TCRM with $\alpha_{\tau}=0.5\tau^{-1}-1$, $0.5\pi\cot(\pi\tau)$ and $-\log(2-2\tau)/\log(1-\tau)$, respectively). The results show that the value of $\omega_{\tau}\xi_{\tau}(\bm{Y}|\bm{x})$  is smaller than that of QR. In addition, $\xi_{\tau}(\bm{Y}|\bm{x})$  with $\text{J}_{\tau}(s)$ in (2.4) can also be used as a new regression model, as shown on the right side of Figure 1c.

\begin{figure}[htbp]
    \centering
    \subfigure[Generalized ES]{
        \includegraphics[width=0.8\textwidth]{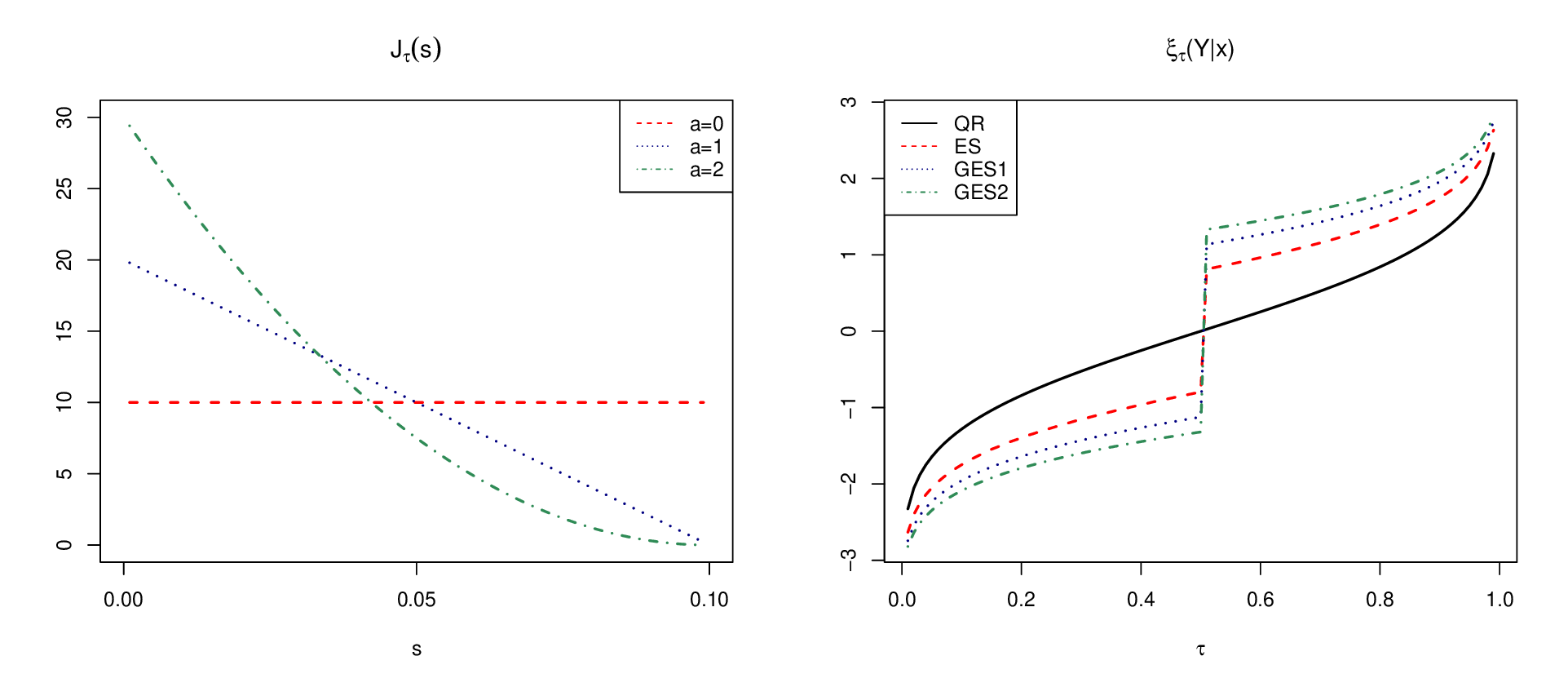}
        \label{fig:sub1}
    }
    \hfill
    \subfigure[Generalized Extremile]{
        \includegraphics[width=0.8\textwidth]{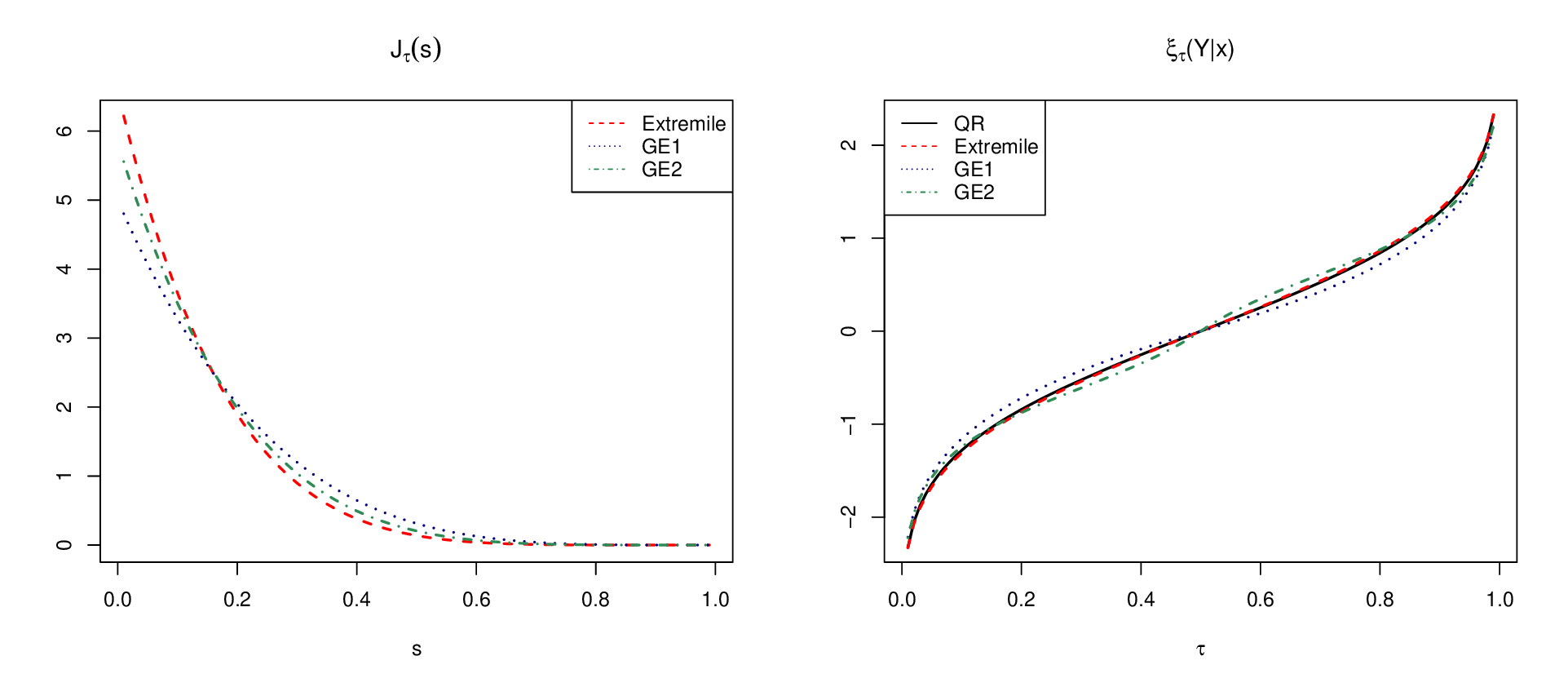}
        \label{fig:sub2}
    }
     \hfill
    \subfigure[Coherent risk measure with truncated Cauchy density function]{
        \includegraphics[width=0.8\textwidth]{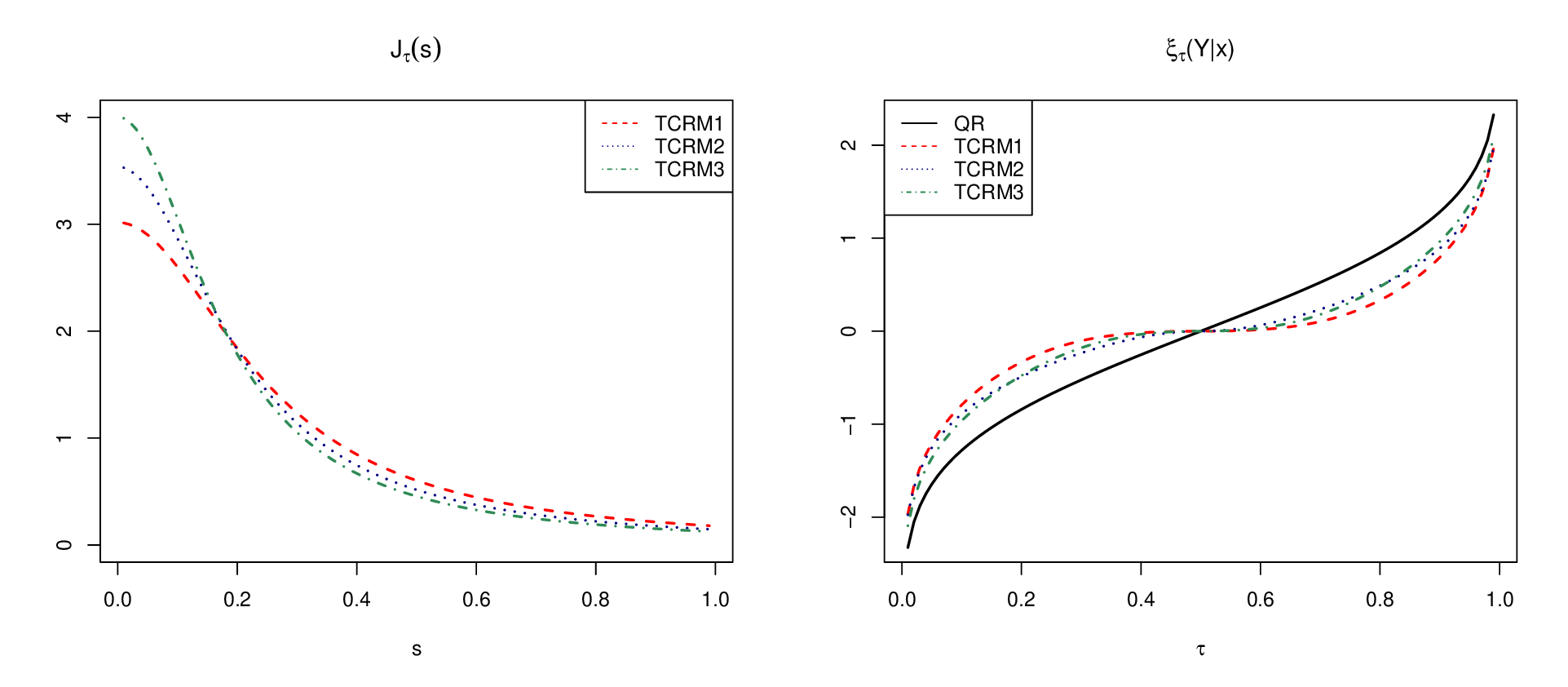}
        \label{fig:sub2}
    }
    \caption{The plots of $\text{J}_{\tau}(s)$s in (2.2)-(2.4) under $\tau=0.1$ and of $\xi_{\tau}(\bm{Y}|\bm{x})$ with $Q_{\bm{Y}|\bm{x}}(\cdot)$ the quantile of the standard normal distribution.   }
    \label{fig:total}
\end{figure}

\begin{table}[htp]
	\footnotesize
	\caption{$\tau$ values corresponding to Copies and Scale for GE and TCRM, respectively.}
	\centering
	\begin{tabular}{@{}cc|ccccc@{}}
\hline
		&Copies&2&5&10&25&50\\
     GE&$\tau$&0.250&0.100&0.050&0.020&0.010\\
		\hline
		&Scale&2&1&1/2&1/4&1/9\\
     TCRM&$\tau$&0.333&0.250&0.167&0.100&0.050\\
	\hline
	\end{tabular}
\end{table}

\begin{figure}[htb]
  \centering
\includegraphics[width=5in]{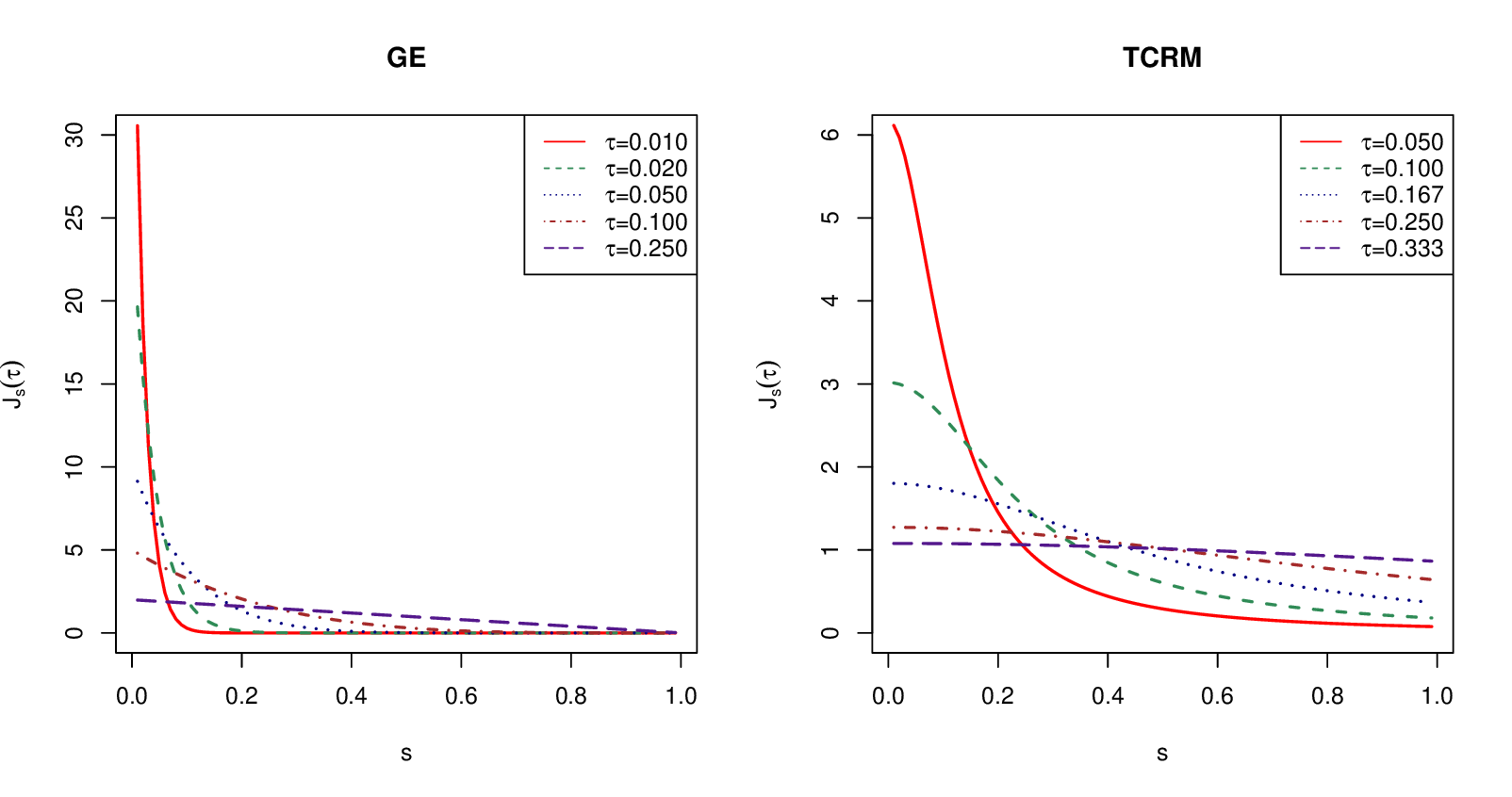}
   \caption{The plots of $\text{J}_{\tau}(s)$ in (2.3) and (2.4) with $\alpha_\tau= 0.5\tau^{-1}-1$ and different $\tau$s according to Table 1. }
\end{figure}
\subsection{Comparing risk measurement tools}
In this section, we explore the relationship between GES, GE, TCRM and QR using the Fr\'{e}chet distribution. The Fr\'{e}chet distribution is one of the extreme value distributions, commonly used for financial risk compared to the Weibull and Gumbel distributions. Risk measurement mainly focuses on the tail situation. In Section 2.2, it focuses on the small quantile situation. Therefore, this section considers another aspect (opportunity), namely the high quantile $\tau\rightarrow 1$.
\begin{thm}
 Suppose the conditions in Theorem 2.1 hold and the conditional distribution of $\bm{Y}$ given $\bm{X}=\bm{x}$ is the Fr\'{e}chet distribution with the distribution function $\exp(-y^{-1/\gamma(\bm{x})})$ on support $[0,+\infty)$ and $\gamma(\bm{x})\in(0,1)$. Then, if also in parts (ii) and (iii), $\lim_{\tau\rightarrow 1}(1-\tau)\alpha_{1-\tau}=A>0$, we have
 \par
 (i) Generalized ES:
 $$\lim_{\tau\rightarrow 1} \frac{\xi_{\tau}(\bm{Y}|\bm{x})}{Q_{\bm{Y}|\bm{x}}(\tau)} =(1+a)B(1-\gamma(\bm{x}),1+a),$$
 \par
 (ii) Generalized Extremile:
 $$\lim_{\tau\rightarrow 1} \frac{\xi_{\tau}(\bm{Y}|\bm{x})}{Q_{\bm{Y}|\bm{x}}(\tau)} =A^{\gamma(\bm{x})}\Gamma(1-\gamma(\bm{x})),$$
  \par
 (iii) TCRM:
 $$\lim_{\tau\rightarrow 1} \frac{\xi_{\tau}(\bm{Y}|\bm{x})}{Q_{\bm{Y}|\bm{x}}(\tau)} =A^{\gamma(\bm{x})}\sec(\gamma(\bm{x})\pi/2),$$
 \\
 where $B(\cdot,\cdot)$ is the beta function and $\Gamma(\cdot)$ is the gamma function.
\end{thm}
\par
The following results can be found from Theorem 2.2:
(1) The larger $\gamma(\cdot)$ is, the thicker the tail will be for Fr\'{e}chet distribution.
From Theorem 2.2(i), ES ($a=0$) is larger than quantile according to $(1+a)B(1-\gamma(\bm{x}),1+a)=\{1-\gamma(\bm{x})\}^{-1}>1$. Moreover, function $(1+a)B(1-\gamma(\bm{x}),1+a)$ increases as $a>0$ increases, so GES1 ($a=1$) and GES2 ($a=2$) behave more conservatively than ES;
(2) For GE and TCRM, the larger A is, the larger GE and TCRM are relative to quantile for $\forall \gamma(\bm{x})\in(0,1)$, that is, the more conservative they are.
For the same $\gamma(\bm{x})\in(0,1)$, GE$>$TCRM according to $\Gamma(1-\gamma(\bm{x}))>\sec(\gamma(\bm{x})\pi/2)$. GE and TCRM are larger than quantile for $A\geq 1$;
(3) Extremiel $A=\ln 2$, and $(\ln 2)^{\gamma(\bm{x})}\Gamma(1-\gamma(\bm{x}))>1$ for $\forall \gamma(\bm{x})\in(0,1)$. Therefore, under the Fr\'{e}chet distribution, Extremiel is more conservative than quantile.
For GE with $\alpha_{\tau}=0.5(1-\tau)^{-1}-1$, then $A=0.5<\ln 2(\approx0.69)$, therefore, Extremiel is also more conservative than this GE. When $\gamma(\bm{x})<0.13168$, GE$<$quantile; otherwise GE$>$quantile;
(4) When $A<1$, there exists a unique $\gamma_0(\bm{x})$ such that when $\gamma(\bm{x})<\gamma_0(\bm{x})$, TCRM$<$quantile; when $\gamma(\bm{x})>\gamma_0(\bm{x})$, TCRM$>$quantile.
Specially, for TCRM with $\alpha_{\tau}=0.5(1-\tau)^{-1}-1$ and $A=0.5$, when $\gamma(\bm{x})<\gamma_0(\bm{x})=0.5$, TCRM$<$quantile; conversely, TCRM$>$quantile.
\par
Next, take specific distributions as examples to compare the commonly used and newly proposed risk measurement tools. We assume that the conditional distribution of $\bm{Y}$ given $\bm{x}$ follows the following six commonly used distributions, which are the t distribution (t(3) and t(1.2), which are in the domain of attraction of the Fr\'{e}chet distribution with 1/3 and 5/6, respectively); standard Normal distribution and exponential distribution (Normal(0,1) and Exp(1), which are in the domain of attraction of the Gumbel distribution); Uniform distribution and Beta distribution (U(0,1) and Beta(2,3), which are in the domain of attraction of the Weibull distribution).
The high quantiles $\tau$ from 0.90 to 0.98 are considered.
We take $a=1$ for GES, $\alpha_{\tau}=0.5\tau^{-1}-1$ for GE and TCRM.
From Figure 3, we can see that:
(1) For all six distributions, the order of their values is
GES$>$ES$>$Extremile$>$GE$>$TCRM, which is consistent with the size of the $\text{J}_{\tau}(s)$ tail weight (Figure 4). Based on this, the appropriate $\text{J}_{\tau}(s)$ can be selected according to risk preference, that is, the larger the $\text{J}_{\tau}(s)$ value at the tail, the greater the risk value (more conservative);
(2) For the distributions in the domain of attraction of the Fr\'{e}chet distribution (t(3) and t(1.2)), the results conform to the conclusion in Theorem 2.2, that is, there is the relationship GES$>$ES$>$Extremile$>$GE$>$TCRM, and Extremile is greater than QR, while TCRM is smaller than QR under t(3) with $\gamma(\bm{x})=1/c<0.5$ and greater than QR under t(1.2) with $\gamma(\bm{x})=5/6>0.5$;
(3) For distributions in the domain of attraction of the Gumbel distribution (Normal(0,1) and Exp(1)), the order of their values is
GES$>$ES$>$Extremile$>$QR$>$GE$>$TCRM, while GES$>$ES$>$QR$>$Extremile$>$GE$>$TCRM
 under distributions in the domain of attraction of the Weibull distribution (U(0,1) and Beta(2,3)). The numerical relationship between Extremile and QR is different here.

 \begin{figure}[htp]
  \centering
\includegraphics[width=5.5in]{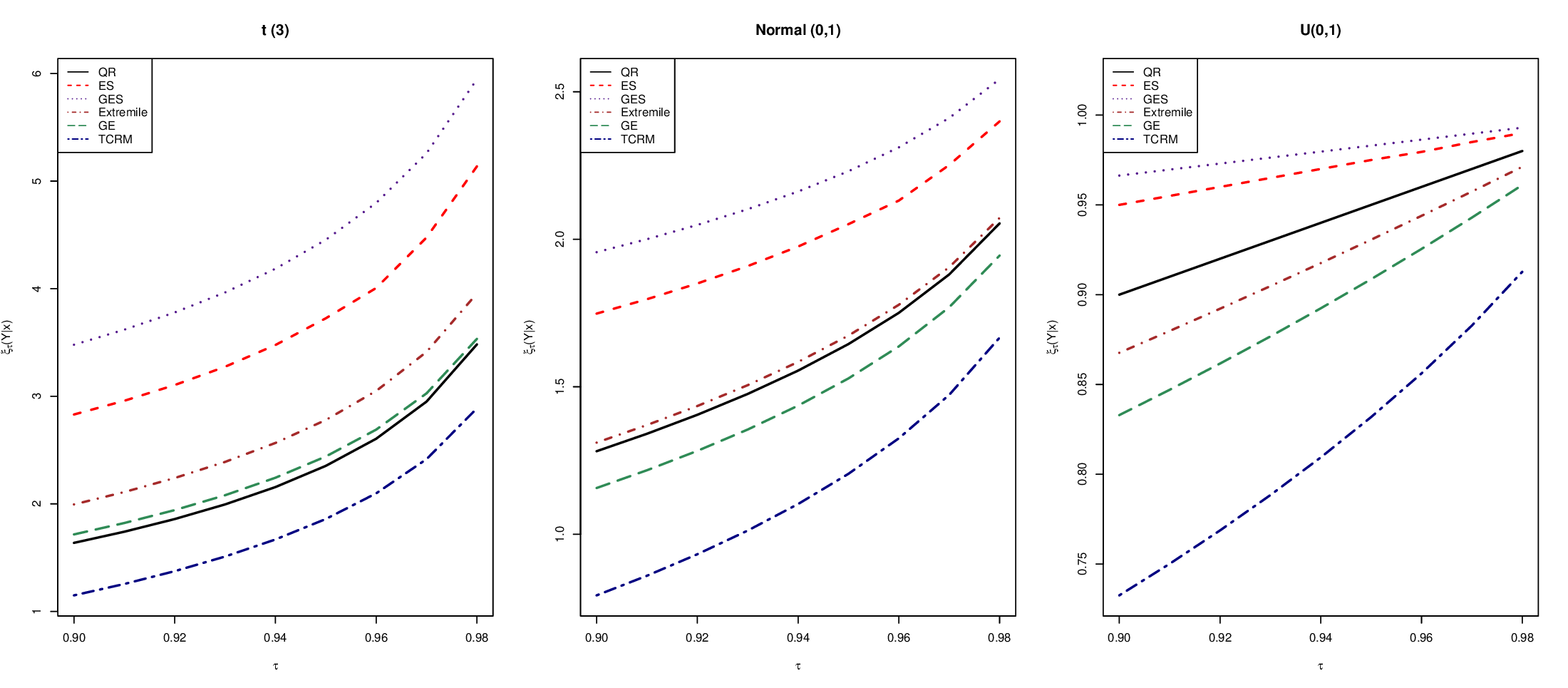}
\includegraphics[width=5.5in]{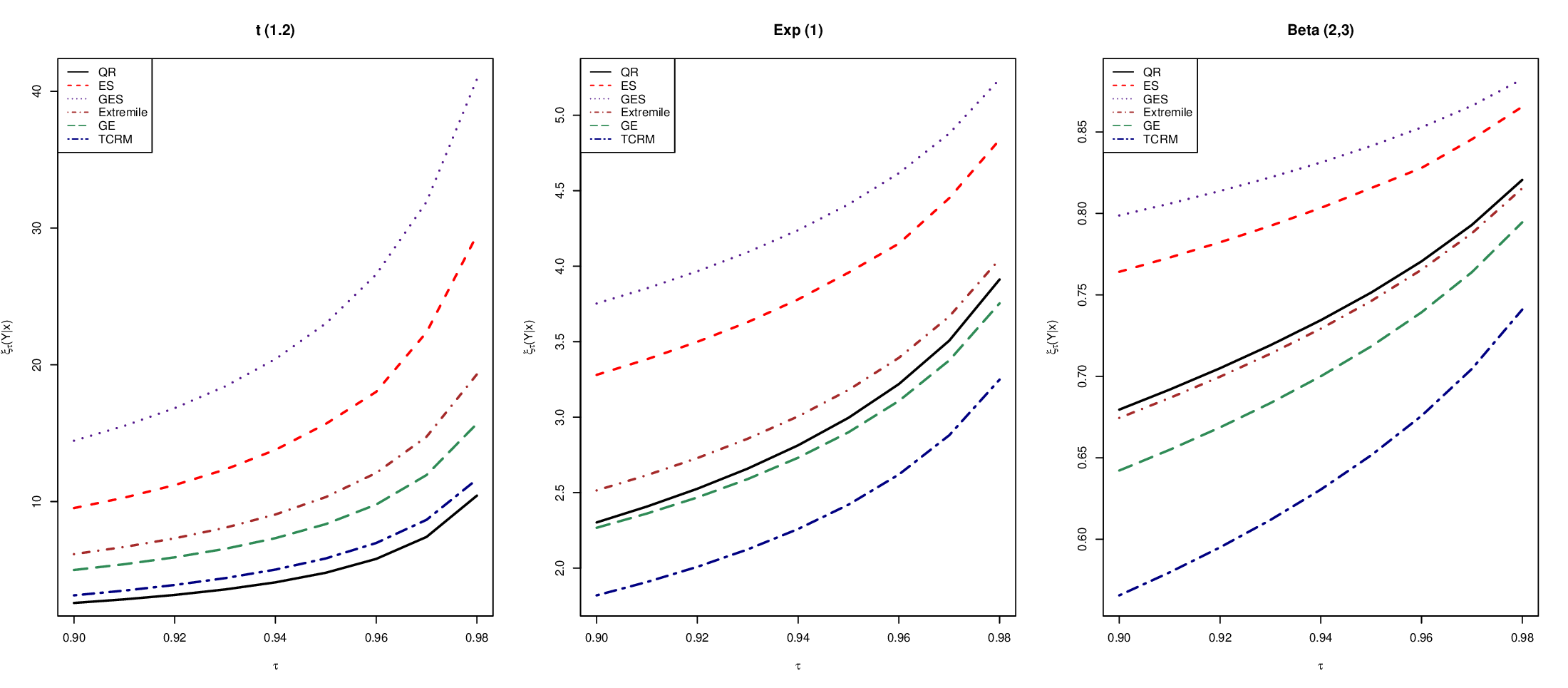}
  \caption{Several risk measures under different commonly used distributions for $\tau\in[0.90,0.98]$.}
\end{figure}

\begin{figure}[htp]
  \centering
 \includegraphics[width=4.6in]{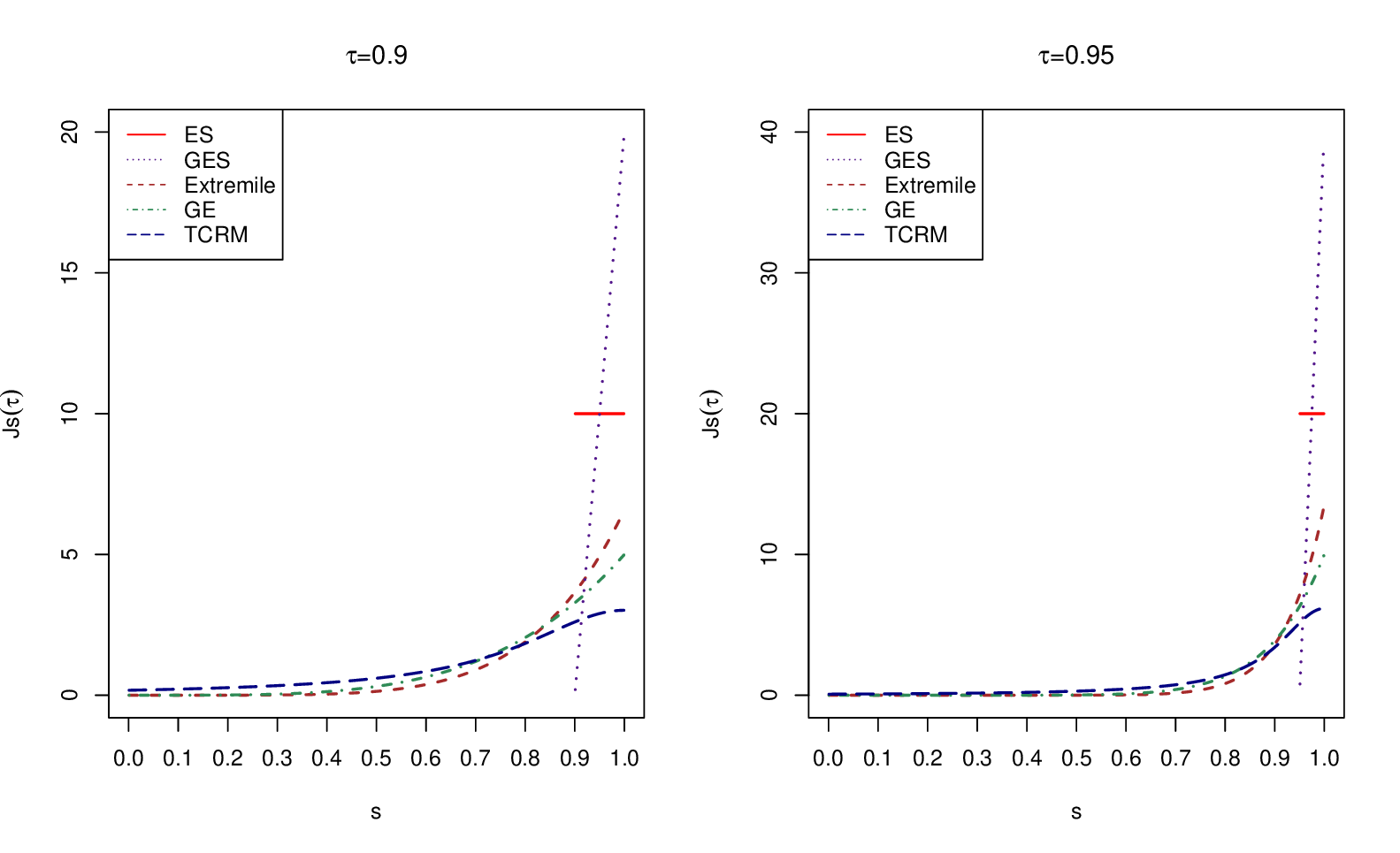}
  \caption{Several weight functions $\text{J}_{\tau}(s)$ under $\tau=0.90,0.95$. }
\end{figure}

\section{Estimation of AQR}
\subsection{Standard estimation method}
In this section, we present a method for estimating $\xi_{\tau}(\bm{Y}|x)$ with univariate $\bm{X}$ for clearly introducing the main idea.
An important approach to the case of multivariate $\bm{X}$ is described in Section 3.2.
Let $\{Y_i, X_i\}_{i=1}^n$ be independent and identically distributed samples from $({\bf Y},{\bf X})$ in model (1.1).
From the last line of equation (2.1), we estimate $\xi_{\tau}(\bm{Y}|x)$ as:
\begin{equation}
	\begin{split}	\hat{\xi}_{\tau}(\bm{Y}|x)=
\int_{0}^{+\infty}\left[1-\text{G}_{\tau}\{\hat{F}(y|x)\}\right]dy
		-\int_{-\infty}^{0}\text{G}_{\tau}\{\hat{F}(y|x)\}dy.
	\end{split}
\end{equation}
In this paper, we use kernel conditional distribution estimation to estimate $F(y|x)$ as:
\begin{equation}
	\begin{split}	\hat{F}(y|x)
={\sum_{i=1}^{n}\text{I}(Y_i\leq y)\text{K}_h\left(X_i-x\right)}/{\sum_{i=1}^{n}
	\text{K}_h\left(X_i-x\right)},
	\end{split}
\end{equation}
where $K_h(\cdot)=K(\cdot/h)/h$, $K(\cdot)$ is a kernel density function and $h>0$ is the bandwidth.
\par
To establish the asymptotic normality of the proposed estimator, the following
technical conditions are imposed.
\par
{\bf C2.} The conditional distribution function $F(y|x)$ has continuous second-order partial derivatives with respect to $x$ and the conditional density function $f(y|x)$ satisfies $c\leq f(y|x)<\infty$ for all $y\in \mathbb{R}$ and $x\in \text{I}_{x}$, where $c$ is a positive constant and $\text{I}_{x}$ is a bounded interval on $\mathbb{R}$.
In addition, the density function $f_{\bm{X}}(\cdot)$ of $\bm{X}$ is positive and continuously differentiable on $\mathbb{R}$.
\par
{\bf C3.} The kernel function $K(\cdot)$ is even, integrable, and twice differentiable with bounded first and second derivatives such that $\int K(u)du=1$, $\int |u^2K(u)|du<\infty$, $\int uK(u)du=0$ and $\int u^2K(u)du\neq0$.
\begin{remark}
Assumption {\bf C2} imposes mild regularity conditions on the conditional distribution and density of $\bm{Y}$ given $\bm{X}$.
Condition {\bf C3} is a mild condition on $K(\cdot)$.  For example, taking a normal density as the kernel function  satisfies condition {\bf C3}.
 \end{remark}

 \begin{thm}
 Assume that $\bm{Y}$ given $\bm{X}=x$ has a finite absolute first moment and that conditions {\bf C1}-{\bf C3} hold. Suppose that $h=n^{-c_1}$ with $c_1\in(1/9,1/5]$ and $n\rightarrow \infty$. Then for given $x\in \text{I}_{x}$, we have
\begin{equation*}
\begin{split}
\sqrt{nh}\left\{\hat{\xi}_{\tau}(\bm{Y}|x)-\xi_{\tau}(\bm{Y}|x)
-\tfrac12 \nu_2^1h^2B_{x}\right\}
\xrightarrow{L}\textrm{N}\left(0,\bm{\Sigma}_{x}\right),
\end{split}
\end{equation*}
where $\nu_a^b=\int_{-\infty}^{+\infty}u^aK^b(u)du$, $B_{x}=-\int_{-\infty}^{+\infty}\text{J}_{\tau}\{F(y|x)\}\{\text{F}''(y|x)
+2\text{F}'(y|x)f'_{\bm{X}}(x)
		/f_{\bm{X}}(x)\}dy,$ $F'(y|x)={\partial F(y|x)}/{\partial x}$,
$F''(y|x)={\partial^2F(y|x)}/{\partial x^2}$, $f'_{\bm{X}}(x)={\partial f_{\bm{X}}(x)}/{\partial x}$, $\xrightarrow{L}$ stands for convergence in distribution, and $$\bm{\Sigma}_{x}=\nu_0^2f^{-1}_{\bm{X}}(x)
\int_{-\infty}^{+\infty}\int_{-\infty}^{+\infty}
\text{J}_{\tau}\{F(y_1|x)\}\text{J}_{\tau}\{F(y_2|x)\}\left\{F(y_1\wedge y_2|x)-F(y_1|x)F(y_2|x)\right\}dy_1dy_2.$$
\end{thm}
\begin{remark}
If we only consider $\bm{Y}$ without the covariable $\bm{X}$, AQR is reduced to $
		\xi_{\tau}(\bm{Y})=\int_0^1Q_{\bm{Y}}(s)\text{J}_{\tau}(s)ds.
$
The corresponding estimate of $\xi_{\tau}(\bm{Y})$ is
$
		\hat{\xi}_{\tau}(\bm{Y})=\int_0^1\hat{Q}_{\bm{Y}}(s)\text{J}_{\tau}(s)ds
=n^{-1}\sum_{i=1}^n\tilde{Y}_{i}\text{J}_{\tau}\{i/(n+1)\},
$
where $\tilde{Y}_{1}\leq\cdots\leq \tilde{Y}_{n}$ denotes the ordered sample and $\hat{Q}_{\bm{Y}}(s)$ is the estimator of $Q_{\bm{Y}}(s)$. For any given $\tau\in(0,1)$ and $\text{E}|\bm{Y}|^{\varsigma}<\infty$ for some $\varsigma>2$, then by Theorem 1(ii) of \cite{r41}, we have
$$
\sqrt{n}\left\{\hat{\xi}_{\tau}(\bm{Y})-\xi_{\tau}(\bm{Y})\right\}
\xrightarrow{L}\textrm{N}\left(0,\int_{0}^{1}\int_{0}^{1}
\text{J}_{\tau}(r)\text{J}_{\tau}(s)(r\wedge s-rs)dF^{-1}(r)dF^{-1}(s)\right).
$$

\end{remark}

\subsection{AQR for distributed systems with large $n$ and $p$}
We first introduce the data in the distributed systems. Define $M={1,\ldots,n}$ as the set of all sample observations where $n$ observations are distributed across $K$ local machines (workers). Decompose $M$ into subsets $M_k$ for $k=1,\ldots,K$, where $M_k$ comprises observations distributed to the $k$th worker, denote $|M_k|=n_k$, and $n$ is the total sample size given by $\sum_{k=1}^Kn_k$.
\par
In the context of distributed systems, the sample size $n$ is often exceedingly large, making it impossible for a single computer to store or run algorithms. Moreover, the data generated in distributed systems is often high-dimensional data.
When the dimensionality of $\bm{X}$ is large the estimation method (3.2) for $F(\cdot|\cdot)$
 will face the ``curse of dimensionality''. Therefore, we assume that there is a $p$-dimensional unknown parameter vector $\bm{\beta}_0$ that makes the following formula true:
\begin{equation}
	\begin{split}
		F(y|\bm{x})=F(y|\bm{x}^{\top}\bm{\beta}_0),
	\end{split}
\end{equation}
where $\bm{x}$ is a p-dimensional vector.
For identification, the
first component of $\bm{\beta}_0$ is positive and $\|\bm{\beta}_0\|_2=1$, where $\|\cdot\|_2$ denotes
the Euclidean 2-norm.
Model (3.3) is the single-index conditional distribution model \citep{r8,r13}.
\par
From the equations (2.1) and (3.3), we can derive $\xi_{\tau}(\bm{Y}|\bm{x})$ as an average
quantile single-index regression (AQSIR) as follows:
\begin{equation}
	\begin{split}
		\xi_{\tau}(\bm{Y}|\bm{x}^{\top}\bm{\beta}_0)		=\int_{0}^{+\infty}\left[1-\text{G}_{\tau}\{F(y|\bm{x}^{\top}\bm{\beta}_0)\}\right]dy
		-\int_{-\infty}^{0}\text{G}_{\tau}\{F(y|\bm{x}^{\top}\bm{\beta}_0)\}dy.
	\end{split}
\end{equation}
According to definitions of (3.3) and (3.4), $\bm{\beta}_0$ is independent of $\tau$. In addition,
 we can estimate $\xi_{\tau}(\bm{Y}|\bm{x}^{\top}\bm{\beta}_0)$ as:
\begin{equation}
	\begin{split}	\hat{\xi}_{\tau}(\bm{Y}|\bm{x}^{\top}\hat{\bm{\beta}})=\int_{0}^{+\infty}
\left[1-\text{G}_{\tau}\{\hat{F}(y|\bm{x}^{\top}\hat{\bm{\beta}})\}\right]dy
		-\int_{-\infty}^{0}\text{G}_{\tau}\{\hat{F}(y|\bm{x}^{\top}\hat{\bm{\beta}})\}dy,
	\end{split}
\end{equation}
where $\hat{\bm{\beta}}$ can be obtained by the
pseudo sum of integrated squares (PSIS) inspired by \cite{r8} and \cite{r9} as:
\begin{equation}
	\begin{split}
		\hat{\bm{\beta}}=&\arg\min_{\bm{\beta}\in\mathbb{R}^p}\frac{1}{n}\sum_{i=1}^n
\int_{-\infty}^{+\infty}\left\{\text{I}(Y_i\leq y)-\hat{F}(y|\bm{X}_i^{\top}\bm{\beta})\right\}^2d\hat{F}(y)\\
=&\arg\min_{\bm{\beta}\in\mathbb{R}^p}\frac{1}{n^2}\sum_{i=1}^n\sum_{j=1}^nL_{ij}(\bm{\beta}),
	\end{split}
\end{equation}
where $\hat{F}(y)=n^{-1}\sum_{i=1}^n\text{I}(Y_i\leq y)$ and $L_{ij}(\bm{\beta})=\{\text{I}(Y_i\leq Y_j)-\hat{F}(Y_j|\bm{X}_i^{\top}\bm{\beta})\}^2$.
\par
Due to the large amount of data $n$ in the distributed system, it is not possible to directly utilize (3.6).
We then use a Taylor expansion of $\bar{L}(\bm{\beta})=n^{-2}
\sum_{i=1}^n\sum_{j=1}^nL_{ij}(\bm{\beta})$ around an initial  estimator $\hat{\bm{\beta}}^{0}$ of $\bm{\beta}_{0}$. This yields:
\begin{equation*}
\begin{split}
\bar{L}(\bm{\beta})=&\bar{L}(\hat{\bm{\beta}}^{0})+(\bm{\beta}-\hat{\bm{\beta}}^{0})^{\top} \nabla\bar{L}(\hat{\bm{\beta}}^{0})+\frac{1}{2}
(\bm{\beta}-\hat{\bm{\beta}}^{0})^{\top}\nabla^2\bar{L}(\hat{\bm{\beta}}^{0})(\bm{\beta}-\hat{\bm{\beta}}^{0})
+o_p(\|\bm{\beta}-\hat{\bm{\beta}}^{0}\|_2^2),
\\
=&\bar{L}(\hat{\bm{\beta}}^{0})+(\bm{\beta}-\hat{\bm{\beta}}^{0})^{\top} \nabla\bar{L}(\hat{\bm{\beta}}^{0})+\frac{1}{2}
(\bm{\beta}-\hat{\bm{\beta}}^{0})^{\top}\nabla^2\bar{L}_1(\hat{\bm{\beta}}^{0})(\bm{\beta}-\hat{\bm{\beta}}^{0})
+o_p(\|\bm{\beta}-\hat{\bm{\beta}}^{0}\|_2^2),
\end{split}
\end{equation*}
where $\nabla$ and $\nabla^2$ represents the first and second derivativeis with respect to $\bm{\beta}$, respectively, $\bar{L}_{1}(\bm{\beta})=n_1^{-2}\sum_{i\in M_1}\sum_{j\in M_1}L_{ij}(\bm{\beta})$ with bandwidth $h_1$ dependent on $n_1$,  and
the last identity in the equation is because $\|\nabla^2\bar{L}(\hat{\bm{\beta}}^{0})-\nabla^2\bar{L}_1(\hat{\bm{\beta}}^{0})\|=o_p(1)$ (see the proof of Theorem 3.2 in the Appendix). The purpose of doing this is to reduce the communication burden, that is, without using the transmission matrix $p\times p$-dimensional matrix $\nabla^2\bar{L}(\hat{\bm{\beta}}^{0})$ ($p\rightarrow\infty$).
\par
Therefore, the estimation of $\bm{\beta}_0$ in model (3.4) can be implemented by solving the following quadratic optimization to obtain something which admits a fast and scalable algorithm to perform optimization under massive and high-dimensional data:
\begin{equation}
	\begin{split}
		\hat{\bm{\beta}}^{1}=&\arg\min_{\bm{\beta}}\left\{\bar{L}(\hat{\bm{\beta}}^{0})+(\bm{\beta}-\hat{\bm{\beta}}^{0})^{\top} \nabla\bar{L}(\hat{\bm{\beta}}^{0})+\frac{1}{2}
(\bm{\beta}-\hat{\bm{\beta}}^{0})^{\top}\nabla^2\bar{L}_1(\hat{\bm{\beta}}^{0})(\bm{\beta}-\hat{\bm{\beta}}^{0})\right\}\\
=&\hat{\bm{\beta}}^{0}-\left\{\nabla^2\bar{L}_{1}(\hat{\bm{\beta}}^{0})\right\}^{-1}
		\nabla\bar{L}(\hat{\bm{\beta}}^{0}),
	\end{split}
\end{equation}
where the initial estimator $\hat{\bm{\beta}}^{0}$ can be obtained by $\hat{\bm{\beta}}^{0}=\arg\min_{\bm{\beta}\in\mathbb{R}^p}\bar{L}_{1}(\bm{\beta})$.
It uses data available only on the first machine, used as the central machine, along with $p$-dimensional gradient vectors $\nabla\bar{L}(\hat{\bm{\beta}}^{0})$ that are sent from the remaining local machines. 
\par
In equation (3.7), $\nabla\bar{L}(\bm{\beta})$ can be split into $n^{-2}\sum_{k_2=1}^K\sum_{i\in M_{k_2}}\sum_{k_1=1}^K\sum_{j\in M_{k_1}}\nabla L_{ij}(\bm{\beta})$ with
$\nabla L_{ij}(\bm{\beta})=-2\{I(Y_i\leq Y_j)-\hat{F}(Y_j|\bm{X}_i^{\top}\bm{\beta})\}\nabla \hat{F}(Y_j|\bm{X}_i^{\top}\bm{\beta})$. Moreover,
$\hat{F}(y|\bm{x}^{\top}{\bm{\beta}})$ and $\nabla \hat{F}(y|\bm{x}^{\top}\bm{\beta})$ are additive, so they can be easily and directly used in distributed systems as:
$	
\hat{F}(y|\bm{x}^{\top}{\bm{\beta}})=S_1(y,\bm{x},\bm{\beta})/S_2(\bm{x},\bm{\beta})
$
and
$
\nabla \hat{F}(y|\bm{x}^{\top}\bm{\beta})=S_3(y,\bm{x},\bm{\beta})/S_2(\bm{x},\bm{\beta})-
 S_1(y,\bm{x},\bm{\beta})\times S_4(\bm{x},\bm{\beta})/S_2^2(\bm{x},\bm{\beta}),
$
where $S_1(y,\bm{x},\bm{\beta})=\sum_{k=1}^K\sum_{i\in M_k}\text{I}(Y_i\leq y)K_h(\bm{X}_i^{\top}\bm{\beta}-\bm{x}^{\top}\bm{\beta})$, $S_2(\bm{x},\bm{\beta})=\linebreak \sum_{k=1}^K\sum_{i\in M_k}K_h(\bm{X}_i^{\top}\bm{\beta}-\bm{x}^{\top}\bm{\beta})$, $S_3(y,\bm{x},\bm{\beta})=\sum_{k=1}^K\sum_{i\in M_k}\text{I}(Y_i\leq y)K'_h(\bm{X}_i^{\top}\bm{\beta}-\bm{x}^{\top}\bm{\beta})(\bm{X}_i-\bm{x})$, and $S_4(\bm{x},\bm{\beta})=\sum_{k=1}^K\sum_{i\in M_k}K'_h(\bm{X}_i^{\top}\bm{\beta}-\bm{x}^{\top}\bm{\beta})(\bm{X}_i-\bm{x})$.
To sum up, (3.7) is communication-efficient. An algorithm for the above distributed estimation method is given in Appendix B of the Supplementary Material.
\par
To establish the asymptotic properties of the proposed estimators, the following
technical conditions are imposed.
\par
{\bf C4.} Suppose that $\inf_{\bm{x}^{\top}\bm{\beta}}f(\bm{x}^{\top}\bm{\beta})>0$ for all $\bm{x}\in \text{I}^p_{\bm{x}}$ and $\bm{\beta}\in \mathbb{R}^p$, where $f(\bm{x}^{\top}\bm{\beta})$ is the density function of $\bm{x}^{\top}\bm{\beta}$.
Moreover, the third derivative of $f(\bm{x}^{\top}\bm{\beta})$ and $\text{E}\{F(y|\bm{X}^{\top}\bm{\beta})(\bm{x}-\bm{X})(\bm{x}-\bm{X})^{\top}|\bm{x}^{\top}\bm{\beta}\}$ with respect to $\bm{x}^{\top}\bm{\beta}$, are Lipschitz continuous in $\bm{x}^{\top}\bm{\beta}$ with the Lipschitz constants being independent of $(y,\bm{x}^{\top}\bm{\beta})$.
\par
{\bf C5.} $\bm{\Sigma}_1=4\text{E}(\bm{A}\bm{A}^{\top})$ with $\bm{A}=\int\left\{\text{I}(Y\leq y)-F(y|\bm{X}^{\top}\bm{\beta}_0)\right\}\nabla F(y|\bm{X}^{\top}\bm{\beta}_0)dF(y)$ and $\bm{\Sigma}_2=2\text{E}\int[\left\{\nabla F(y|\bm{X}^{\top}\bm{\beta}_0)\right\}^2-\{\text{I}(Y\leq y)-F(y|\bm{X}^{\top}\bm{\beta}_0)\}\nabla^2 F(y|\bm{X}^{\top}\bm{\beta}_0)]dF(y)$ are non-singular. Moreover, the minimum eigenvalue of $\bm{\Sigma}_2$ is positive.
\begin{remark}
The condition {\bf C4} is the smoothness condition required for the uniqueness and 
convergence of the estimator.
{\bf C5} is to ensure the asymptotic normality of the estimator.
 Conditions {\bf C4} and {\bf C5} are the general conditions for establishing the 
 consistency and asymptotic normality of the single-index conditional distribution model (3.3) \citep{r8,r13}.
\end{remark}

\begin{thm}
	Assume that $\bm{Y}$ given $\bm{X}=\bm{x}$ has a finite absolute first moment and that conditions {\bf C1}-{\bf C5} are satisfied.
	Suppose we have an initial estimator $\hat{\bm{\beta}}^{0}$ with
	$\|\hat{\bm{\beta}}^{0}-\bm{\beta}_0\|_2=O_p(n_1^{-1/2})$, $h=O(n^{-c_2})$ and $h_1=O(n_1^{-c_2})$ with $c_2\in (1/8,1/5)$. Then, under $n_1\rightarrow\infty$, we have
	\begin{equation*}
		\begin{split}
			\|\hat{\bm{\beta}}^{1}-\bm{\beta}_0\|_2
			=O_p(n^{-1/2})+o_p\left(n_1^{-1/2}\cdot\sqrt{\frac{\ln n_1}{n_1h_1^5}}\right).
		\end{split}
	\end{equation*}
Furthermore, for the multiple rounds estimator $\hat{\bm{\beta}}^{q}
=\hat{\bm{\beta}}^{q-1}-\{\nabla^2\bar{L}_{1}(\hat{\bm{\beta}}^{q-1})\}^{-1}
		\nabla\bar{L}(\hat{\bm{\beta}}^{q-1})$ with
$q\geq[\ln (n/n_1)/\ln (n_1h_1^5/\ln n_1)]$, we have
	\begin{equation*}
		\begin{split}
		&(i)~\|\hat{\bm{\beta}}^{q}-\bm{\beta}_0\|_2=O_p(n^{-1/2}),\\
&(ii)~\|\hat{\bm{\beta}}^{q}-\hat{\bm{\beta}}\|_2=O_p\left(n^{-1/2}\cdot\sqrt{\frac{\ln n_1}{n_1h_1^5}}\right)+o_p\left(n_1^{-1/2}\cdot\left\{\frac{\ln n_1}{n_1h_1^5}\right\}^{q/2}\right)=o_p(n^{-1/2}),\\
&(iii)~\sqrt{n}(\hat{\bm{\beta}}^{q}-\bm{\beta}_0)
	\xrightarrow{L}\mathcal{N}(\bm{0},\bm{\Sigma}_2^{-1}\bm{\Sigma}_1\bm{\Sigma}_2^{-1}).
	\end{split}
	\end{equation*}
\end{thm}
  \begin{thm}
Suppose the conditions in Theorem 3.2 hold. Then, we have
\begin{equation*}
\begin{split}
\sqrt{nh}\left\{\hat{\xi}_{\tau}(\bm{Y}|\bm{x}^{\top}\hat{\bm{\beta}}^q)-
\xi_{\tau}(\bm{Y}|\bm{x}^{\top}\bm{\beta}_0)
-\tfrac12 \nu_2^1h^2B_{\bm{x}^{\top}\bm{\beta}_0}\right\}
\xrightarrow{L}\textrm{N}\left(0,\bm{\Sigma}_{\bm{x}^{\top}\bm{\beta}_0}\right),
\end{split}
\end{equation*}
where $B_{\bm{x}^{\top}\bm{\beta}_0}=-\int_{-\infty}^{+\infty}\text{J}_{\tau}\{F(y|\bm{x}^{\top}\bm{\beta}_0)\}
\left\{F''(y|\bm{x}^{\top}\bm{\beta}_0)+2F'(y|\bm{x}^{\top}\bm{\beta}_0)\frac{f'_{\bm{X}^{\top}\bm{\beta}_0}
(\bm{x}^{\top}\bm{\beta}_0)}{f_{\bm{X}^{\top}\bm{\beta}_0}(\bm{x}^{\top}\bm{\beta}_0)}\right\}dy$ and \begin{equation*}
\begin{split}
\bm{\Sigma}_{\bm{x}^{\top}\bm{\beta}_0}=&\nu_0^2f^{-1}_{\bm{X}^{\top}
\bm{\beta}_0}(\bm{x}^{\top}\bm{\beta}_0)\int_{-\infty}^{+\infty}\int_{-\infty}^{+\infty}
\text{J}_{\tau}\{F(y_1|\bm{x}^{\top}\bm{\beta}_0)\}\text{J}_{\tau}\{F(y_2|\bm{x}^{\top}\bm{\beta}_0)\}\\
&\times
\left\{F(y_1\wedge y_2|\bm{x}^{\top}\bm{\beta}_0)-F(y_1|\bm{x}^{\top}\bm{\beta}_0)F(y_2|\bm{x}^{\top}\bm{\beta}_0)\right\}dy_1dy_2.
\end{split}
\end{equation*}

\end{thm}

\section{Numerical studies}
In this section, we first employ Monte Carlo simulation studies to evaluate the finite-sample performance of the proposed procedures. Subsequently, we illustrate the application of the proposed methods through two real-data analyses. The versions of AQR considered here are identical to those described in Section 2.3. The standard normal density is utilized as the kernel function, and the bandwidth $h$ is determined via the cross-validation method \citep{r38} in this section. All programs are implemented using \textsf{R} code.

\subsection{Simulation example 1: standard estimation method}
In this subsection, we study the estimation method proposed in Section 3.1 for the regression models and risk measures (Section 2) involved in AQR.
We generate 300 data points from the model:
$
\bm{Y}=20\sin(\pi\bm{X})+\bm{\varepsilon},
$
where $\bm{X}$ is drawn from a normal distribution $N(0,1)$.
Three error distributions of $\bm{\varepsilon}$ are considered:
Normal(0,1), t(3) and Exp(1), where Normal(0,1) is the most commonly used, t(3) is a thick-tailed distribution and Exp(1) is one-sided. The versions of AQR considered are the same as in Section 2.3, including choices of $a$ and $\alpha_\tau$.
\par
The relative percentage absolute deviation (RPAD) is used to assess the performance of estimates as:
$
\text{RPAD}={|\hat{\xi}_{\tau}(\bm{Y}|x)-\xi_{\tau}(\bm{Y}|x)|}/{|\xi_{\tau}(\bm{Y}|x)|}\times 100\%.
$
We take $x=-0.5$ for small values $\tau=0.05,0.10$ and $x=0.5$ for large values  $\tau=0.90,0.95$, respectively. They represent loss (negative) and gain (positive). Simulation results are all the average of 500 simulation replications.
 For extreme values $\tau\in\{0.05,0.10,0.90,0.95\}$ in Table 2, since all \text{RPAD} values are less than $10\%$ (most are less than $5\%$), the proposed estimation method performs well.

\begin{table}[htp]
	\footnotesize
	\caption{The mean and standard deviation (in parentheses) of 500 replicates of RPADs (\%) for $\tau=0.05,0.10,0.90,0.95$ under error following Normal(0,1), t(3) and Exp(1).}
	\centering
	\begin{tabular}{@{}cl|cccc@{}}
		\hline
		Error&Method&$\tau$=0.05&$\tau$=0.1&$\tau$=0.9&$\tau$=0.95\\
		\hline
  Normal(0,1)&        ES&3.14~(2.08)&2.61~(1.86)&2.73~(1.80)&3.36~(2.03)\\
         &GES&3.59~(2.20)&2.97~(2.01)&3.17~(1.95)&3.90~(2.14)\\
   &Extremile&2.54~(1.76)&2.17~(1.53)&2.19~(1.52)&2.66~(1.73)\\
          &GE&2.36~(1.66)&2.08~(1.45)&2.06~(1.45)&2.45~(1.63)\\
        &TCRM&2.15~(1.15)&1.95~(1.39)&1.90~(1.34)&2.18~(1.49)\\
		\hline
      t(3)&    ES&6.99~(4.63)&5.05~(3.72)&5.73~(3.97)&7.74~(5.25)\\
        & GES&8.77~(5.68)&6.40~(4.35)&7.15~(4.71)&9.33~(5.43)\\
  & Extremile&4.97~(3.46)&3.65~(2.57)&4.04~(2.62)&5.55~(3.63)\\
          &GE&4.33~(3.07)&3.27~(2.28)&3.57~(2.31)&4.83~(3.19)\\
       & TCRM&3.59~(2.52)&2.78~(1.94)&2.97~(1.97)&3.96~(2.60)\\
		\hline
    Exp(1)&      ES&0.94~(0.78)&0.81~(0.99)&4.43~(2.84)&5.64~(3.22)\\
        & GES&0.95~(0.76)&0.96~(0.78)&5.25~(3.07)&6.72~(3.50)\\
   &Extremile&1.03~(0.81)&1.18~(0.90)&3.27~(2.19)&4.30~(2.72)\\
          &GE&1.07~(0.83)&1.25~(0.96)&2.94~(1.99)&3.86~(2.52)\\
        &TCRM&1.19~(0.92)&1.45~(1.19)&2.44~(1.65)&3.21~(2.16)\\
	\hline
	\end{tabular}
\end{table}

\subsection{Simulation example 2: distributed estimation method}
In this subsection, we study the distributed estimation method proposed in Section 3.2 for the regression models and risk measures involved in AQR.
We generate sample data from the model:
$
\bm{Y}=(\bm{X}^{\top}\bm{\beta}_0)^2+\bm{\varepsilon},
$
where $\bm{X}=(\bm{X}_1,\bm{X}_2)^{\top}$ are drawn from a normal distribution Normal(2,1),  $\bm{\beta}_0=(1,2)^{\top}/\sqrt{5}$ and $\bm{\varepsilon}$ follows a standard normal distribution.
For the distributed estimation method, we set the number of machines to 10 with a sample size of 50 on each machine. In addition, we set ALL to be the estimator directly using all 500 data points. We take
 $\bm{x}=(2,2)^{\top}$ and values of $\tau=0.1,0.9$ are considered. Simulation results are all the average of 100 simulation replications.
\par
We evaluated parameter estimation performance using the average absolute error (AAE) criterion: $\text{AAE}=\sum_{j=1}^2|\bm{\hat{\beta}}_j-\bm{\beta}_{0,j}|/2$, where $\bm{\hat{\beta}}$ is obtained by (3.6) for the all-data learning (ALL) method and equation (3.7) for the distributed estimation (DE) approach.
Analysis of 100 simulation replicates revealed comparable performance between methods: ALL achieved a mean AAE of 0.0178 (standard deviation=0.0435) while DE showed a mean AAE of 0.0381 (standard deviation=0.0504). Both estimation strategies demonstrated robust performance across parameter configurations.
For scenarios with $\tau=0.1$ and $\tau=0.9$ (Table 3), all RPAD values remained at 9\% or below, indicating satisfactory estimation accuracy. Notably, the distributed method's performance metrics (AAE and RPAD) closely approximate those of the direct all-data approach. These results confirm the proposed distributed estimation framework maintains statistical efficacy while offering computational advantages through data partitioning.

\begin{table}[htp]
	\footnotesize
	\caption{The mean and standard deviation (in parentheses) of 100 replicates of RPADs (\%) for different methods with $\tau=0.1$ and 0.9.}
	\centering
	\begin{tabular}{@{}l|cc|cc@{}}
		\hline
&$\tau=0.1$&&$\tau=0.9$\\
\cline{2-5}
		Method&ALL&DE&ALL&DE\\
		\hline
       ES&4.85~(4.56)&7.83~(9.86)&3.69~(3.15)&6.17~(8.31)\\
      GES&5.48~(5.32)&9.00~(11.4)&3.87~(3.38)&6.60~(9.29)\\
Extremile&3.60~(3.07)&5.50~(6.65)&3.06~(2.48)&5.11~(6.53)\\
       GE&3.24~(2.65)&4.82~(5.68)&2.84~(2.31)&4.71~(5.90)\\
     TCRM&2.49~(1.88)&3.40~(3.51)&2.36~(1.94)&3.80~(4.46)\\
	\hline
	\end{tabular}
\end{table}

\subsection{Real data example 1: Investment portfolio}
In this section, AQR is applied to investment portfolios to illustrate its practical application in the financial field.
The 10 stocks in the portfolio, with reference to the Blackrock U.S. Flexible Equity Fund (BR), are MSFT, AMZN, META, V, NVDA, CIEN, ICE, APD, CAH, WFC.
The 250-trading day data set in 2023 is used as the fit set, while the 252-trading day data set in 2024 is used as the test set. The 10 stock data is downloaded
from the website of Yahoo Finance (https://hk.finance.yahoo.com). The BR has performed well with returns of 22.46\% and 15.21\% in 2023 and 2024, respectively.
For specific information about BR see
https://www.blackrock.com/cn/products/228610/bgf-us-flexible-equity-fund-a2-usd.
\par
The weight $\bm{\alpha}$ of the specific portfolio is chosen to minimize the AQR for a specified value of $\tau$:
\begin{equation*}
\begin{split}
\min_{\bm{\alpha}}~\omega_{\tau}\xi_{\tau}(\bm{\alpha}^{\top}\bm{Y}),~~~
s.t.~\bm{\alpha}^{\top}\bm{1}=1, \bm{\alpha}\geq 0,
\end{split}
\end{equation*}
where $\xi_{\tau}(\bm{\alpha}\bm{Y})$ is defined and estimated in Remark 3.2, $\bm{1}$ is a $10\times 1$ dimensional vector with all 1 elements, and $\bm{Y}=(Y_1,\ldots,Y_{10})^{\top}$ is the logarithmic return of the above 10 stocks. The optimal $\bm{\alpha}$ under $\tau=0.05$ for different risks (mentioned in Section 2) is calculated, as shown in Table 4. BR results are excluded from Table 4 due to its weighting parameter not equaling 1.  This limitation arises because the fund company discloses only its top 10 holdings, withholding all additional portfolio information.
\par
Portfolio performance is evaluated using two complementary metrics: (1) the Sharpe Ratio (SR), calculated as annualized return divided by return volatility (standard deviation), and (2) the Percentage of Days (PD) with excess returns relative to the benchmark (BR). Both indicators follow a ``higher-is-better'' paradigm. Analysis of Table 5 reveals that the TCRM strategy demonstrates superior performance across both evaluation dimensions, achieving top-ranked SR and PD values among all seven methodologies examined. Notably, TCRM generates a 52.27\% absolute return, outperforming BR's 15.21\% by a margin of 3706 basis points. These results collectively validate TCRM as the optimal choice under the specified evaluation framework.

\begin{table}[htp]
	\footnotesize
	\caption{The values of $\bm{\alpha}$ based on fit data for different methods with $\tau=0.05$.}
	\centering
	\begin{tabular}{@{}l|cccccccccc@{}}
		\hline
    Method&MSFT&AMZN&META&V&NVDA&CIEN&ICE&APD&CAH&WFC\\
		\hline
       QR&0.140& 0.014& 0.015& 0.088& 0.230& 0.066&0.251& 0.110& 0.082&0.004\\
       ES&0.072& 0.018& 0.000& 0.140& 0.020& 0.100& 0.190& 0.190& 0.063&0.207\\
      GES&0.244& 0.120& 0.042& 0.110& 0.110& 0.033& 0.150& 0.031& 0.150&0.010\\
Extremile&0.006& 0.058& 0.099& 0.120& 0.011& 0.110& 0.090& 0.018&0.368& 0.120\\
       GE&0.068& 0.110& 0.030& 0.137& 0.063&0.022& 0.130&0.150& 0.140& 0.150\\
     TCRM&0.061& 0.050& 0.190& 0.021& 0.140& 0.047& 0.100& 0.028& 0.100& 0.263\\
	\hline
	\end{tabular}
\end{table}

\begin{table}[htp]
	\footnotesize
	\caption{The SRs and PDs (\%) based on test data for different methods with $\tau=0.05$.}
	\centering
	\begin{tabular}{@{}c|cccccccccc@{}}
		\hline
Method&BR&QR&ES&GES&Extremile&GE&TCRM\\
		\hline
SR&18.03&41.99&33.63&36.41&39.32&38.41&47.31\\
PD&-&53.17&49.21&51.98&50.79&50.40&56.75\\
		\hline
	\end{tabular}
\end{table}

\subsection{Real data example 2: Beijing multi-site air quality dataset}
We apply the proposed AQSIR and its distributed  estimator method in Section 3.2 to the analysis of a Beijing multi-site air quality dataset \citep{r33}.
This dataset includes air pollutant (PM$_{2.5}$) data from 12 nationally controlled air quality monitoring sites. The air quality data are from the Beijing Municipal Environmental Monitoring Center. The standard value of PM$_{2.5}$ in China is that the 24-hour average concentration is less than 75 micrograms per cubic meter.
It is stipulated that a 24-hour average of up to 35 micrograms per cubic meter  is optimal, and up to 75 micrograms per cubic meter is good. If it exceeds 75, it constitutes pollution. The dataset can be obtained from https://archive.ics.uci.edu/dataset/501/beijing+multi+site+air+quality+data.
\par
The official air quality statistics in China are predicated on daily PM$_{2.5}$ values. Nevertheless, it is known that the observed PM$_{2.5}$ levels are affected by meteorological conditions \citep{r34}. Secondary generation of fine PMs is an act of interaction, such as being significantly promoted by high humidity combined with high temperature and calm wind.
Therefore, the meteorological data for each air quality site are obtained from the nearest weather station of the China Meteorological Administration.
The meteorological variables are  temperature (TEMP), pressure (PRES), dew point temperature (DEWP) and wind speed (WSPM).
\cite{r34} used a non-parametric mean regression model to analyze the dataset, which is a special case of AQSIR with $\tau=0.5$.
\par
The histogram of PM$_{2.5}$ in Figure 5 reveals a significant right skew. Moreover, people tend to be more concerned about high PM$_{2.5}$ levels rather than the average. Therefore, it is more appropriate to analyze this dataset using a non-mean regression model such as AQSIR. This section focuses on daily data from the winter of 2016/17 (December 2016 to February 2017), as winter typically exhibits the highest average PM$_{2.5}$ levels compared to other seasons. We therefore have $K=12$, $n_1=\cdots=n_{12}=90$ and $n=1080.$
In addition, in order to eliminate the differences in the range of variation of covariates, data standardization was performed on the four meteorological variables.

\begin{figure}[htp]
  \centering
  \includegraphics[width=3in]{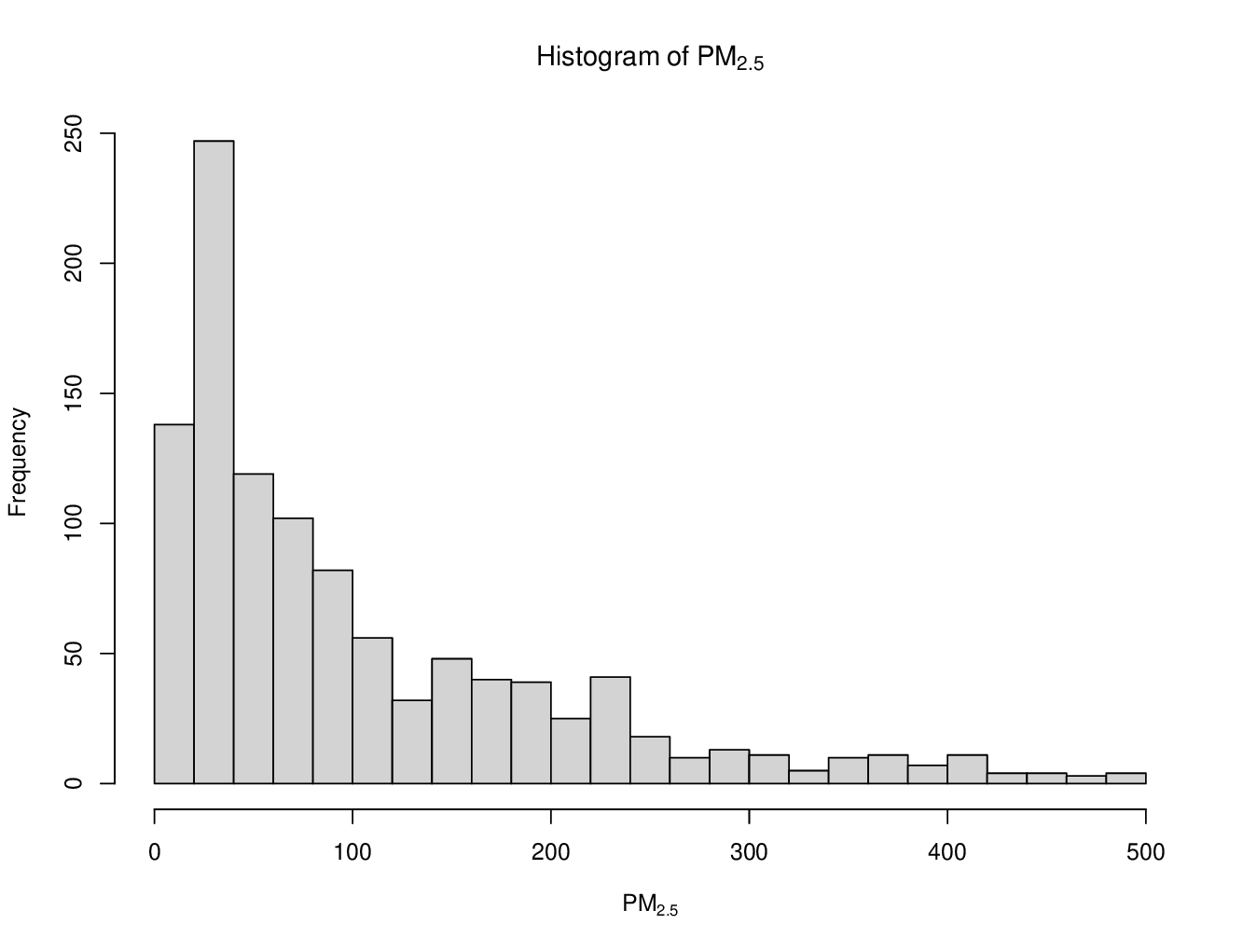}
  \caption{Histogram of PM$_{2.5}$ in the Beijing multi-site air quality dataset.}
  \end{figure}

\par
First, however, we use the proposed AQSIR to analyze the full dataset ignoring its distributed structure. Inspired by the functions (4.2) and (4.3) in \cite{r34}, we calculate the
$
\text{Average PM}_{2.5}=1080^{-1}\sum_{i=1}^{1080}\hat{\xi}_{\tau}(\bm{Y}|\bm{X}_i^{\top}\hat{\bm{\beta}}),
$
where $\bm{Y}$ is $\text{PM}_{2.5}$, $\bm{X}$ is (TEMP, PRES, DEWP, WSPM),
$\hat{\xi}_{\tau}(\bm{Y}|\cdot)$ is obtained by (3.5)
and $\hat{\bm{\beta}}=(0.370,0.275,-0.814,0.354)^{\top}$ is obtained by (3.6). The Average PM$_{2.5}$s under $\tau$ from 0.01 to 0.99 are presented in Table 6. The regularities shown by the results are consistent with those analyzed in Section 2. Moreover, from Table 6, it can be seen that  the median (83, QR) and mean (98, \cite{r34}) of PM$_{2.5}$ in Beijing are greater than the standard value 75.

\begin{table}[htp]
	\footnotesize
\caption{Average PM$_{2.5}$ with $\tau=0.01$ to 0.99 for Beijing multi-site air quality dataset.}
	\centering
	\begin{tabular}{@{}l|ccccccccccccc@{}}
		\hline
Method&$\tau=0.01$&0.05&0.1&0.2&0.3&0.4&0.5&0.6&0.7&0.8&0.9&0.95&0.99\\
\hline
QR       &9&19&27&41&55&68&83&99&117&145&193&230&309 \\
Extremile&10&21&30&47&63&80&98&113&132&158&199&238&311\\
GE       &11&24&36&54&70&84&98&109&124&146&184&222&298\\
TCRM     &17&37&53&77&91&97&98&99&104&119&155&192&272\\
	\hline
	\end{tabular}
\end{table}

\begin{table}[htp]
	\footnotesize
	\caption{ Intervals of $\tau$ (from 0.01 to 0.99) corresponding to Average PM$_{2.5}$ for different methods with Beijing multi-site air quality dataset, where 499 is the maximum.}
	\centering
	\begin{tabular}{@{}c|cccccccc@{}}
		\hline
		Average PM$_{2.5}$&QR&ES&GES&Extremile&GE&TCRM\\
   \hline
(0,35)&  (0.01,0.16)&  (0.01,0.32)&  (0.01,0.49)&(0.01,0.13)&  (0.01,0.10)&(0.01,0.05) \\
(35,75)& (0.16,0.45)&   (0.32,0.50)&   (0.49,0.50)& (0.13,0.37)& (0.10,0.34)&(0.05,0.19)\\
(75,499)&(0.45,0.99)&   (0.50,0.99)&   (0.50,0.99)& (0.37,0.99)&(0.34,0.99)&(0.19,0.99)\\
 \hline
	\end{tabular}
\end{table}
\par
In order to extract more information from the data, we use non-mean regression models like AQSIR
 for further analysis. As can be seen from Table 7, (i) the results based on GES are not very reasonable, because the pollution is most severe in winter, but with GES as the standard, there is too much optimal and almost no good; (ii) TCRM, on the other hand, have too few optimal and good days;
(iii) QR and ES results are more optimistic, with about half optimal and good days. However, the year-round good days are about 50\%, and winter should be a little lower;
(iv) The results for Extremile and GE are close and reasonable.
The essence of both methods is the same, the 35 cut-off point is about $\text{E}\left\{\min(\bm{Y}_{\bm{x}}^1,\ldots,\bm{Y}_{\bm{x}}^{5})\right\}$, and the 75 cut-off point is
between $\text{E}(\bm{Y}|\bm{x})$ and $\text{E}\left\{\min(\bm{Y}_{\bm{x}}^1,\bm{Y}_{\bm{x}}^{2})\right\}$.
Based on the above analysis, GE is a good choice for Beijing multi-site air quality dataset. Because Extremile is a special case, GE explains $\tau$ more visually than Extremile; (v) finally, taking GE as an indicator, under $\tau=0.1$, Average PM$_{2.5}$=36 indicates that there are fewer optimal days in winter. In contrast, under $\tau=0.9$, Average PM$_{2.5}$=184 is about 2.5 times the critical pollution value of 75, which belongs to severe pollution, and even under $\tau=0.99$, serious pollution (Average PM$_{2.5}$=298) is reached, although it rarely occurs.
Such serious pollution has aroused the attention of the relevant departments who have attempted to take control of PM$_{2.5}$ pollution.
In the most recent winter (December 2023 to February 2024), the average concentration of PM$_{2.5}$=38 was close to the optimal threshold of 35, and much smaller than 98 (the winter of 2016/17).
Moreover, the average annual concentration of PM$_{2.5}$ in Beijing's atmospheric environment in 2023 was 32 micrograms per cubic meter, and the average annual concentration of PM${2.5}$ in Beijing in the first three quarters of 2024 (January-September) was 29.
The above data is from the Beijing Municipal Ecology and Environment Bureau (https://sthjj.beijing.gov.cn/bjhrb/index/index.html).
\par
Finally, we consider the distributed method in Section 3.2, because the data comes from 12 nationally controlled air quality monitoring sites. The initial estimator in our method is based on the first site (Aotizhongxin) and the number of ``machines'' is $K=12$. Table 8 lists the absolute deviation between Average PM$_{2.5}$ based on the all data analysis method and Average PM$_{2.5}$ based on the distributed method under different $\tau$s. It can be seen that the distributed method gives results very close to those of the all data analysis method because the absolute deviation between the two is small compared to the values of the full data analysis in Table 6. Therefore, the proposed distributed method is effective.

\begin{table}[htp]
	\footnotesize
	\caption{Absolute deviation between all data and distributed methods estimates of Average PM$_{2.5}$ with $\tau=0.01$ to 0.99 for Beijing multi-site air quality dataset. }
	\centering
	\begin{tabular}{@{}l|ccccccccccccc@{}}
		\hline
Method&$\tau=0.01$&0.05&0.1&0.2&0.3&0.4&0.5&0.6&0.7&0.8&0.9&0.95&0.99\\
\hline
QR&0.15&0.02&0.14&0.43&0.42&0.54&0.46&0.13&0.14&0.16&0.06&0.26&0.89 \\
Extremile&0.08&0.07&0.19&0.27&0.27&0.25&0.20&0.17&0.12&0.03&0.15& 0.38&1.19\\
GE&0.07&0.12&0.23&0.28&0.27&0.24&0.20&0.18&0.15&0.07&0.08& 0.28& 0.99\\
TCRM&0.01&0.17&  0.24&  0.24&  0.21&  0.20&  0.20&  0.19&  0.18&  0.15&  0.03&  0.14&  0.75\\
	\hline
	\end{tabular}
\end{table}

\section{Conclusion}
The article introduces a novel family of non-mean regression models, termed Average
 Quantile Regression (AQR), which also functions as a coherent risk measure through an appropriately defined averaging function, $\text{J}_{\tau}(s)$. Although certain conditions are imposed on $\text{J}_{\tau}(s)$, the flexibility and adaptability of this function, depending on both $\tau$ and $s$, enable AQR to encompass many classical and recently proposed regression models and risk measures as special cases.
\par
Section 2 presents several new non-mean regression models and coherent risk measures, offering a foundational framework for readers to develop their own models and risk tools based on AQR. Examples demonstrate the applicability of AQR to the analysis of high-dimensional and large-scale datasets, particularly those generated by distributed systems. These applications highlight AQR's potential for extension to other complex data types, such as streaming data.
Streaming data, a prominent form of big data, is characterized by continuously arriving, sequentially dependent observations that accumulate over time. Analyzing such data requires updatable and memory-efficient processing methods \citep{r32}. To apply AQR in streaming settings, a local polynomial interpolation method \citep{r29} can be employed to obtain the online update estimator $\hat{F}(y|\bm{x})$ in Equation (3.2), which in turn allows for real-time updating of the quantile estimator $\hat{\xi}_{\tau}(\bm{Y}|\bm{x})$ as new data become available.

\section*{Acknowledgments}
This research is supported by the Ministry of Education of the People's Republic of China, Humanities and Social Science Foundation (Grant No.22YJC910005), National Social Science Foundation of China (Grant No.20BTJ049, No.21BTU040), National Natural Science Foundation of China (Grant No.U23A2064), and Zhejiang Provincial Natural Science Foundation (Grant No.LY24A010004).

\section*{Supplementary material}
The proofs of theorems and algorithm are given in the Supplementary Material file.

\setlength{\bibsep}{0pt}
\bibliographystyle{Chicago}
\bibliography{ref}
\end{document}